\documentclass[sn-mathphys-num]{sn-jnl}

\usepackage{graphicx}%
\usepackage{multirow}%
\usepackage{amsmath,amssymb,amsfonts}%
\usepackage{bbm}
\usepackage{amsthm}%
\usepackage{mathrsfs}%
\usepackage[title]{appendix}%
\usepackage{xcolor}%
\usepackage{textcomp}%
\usepackage{manyfoot}%
\usepackage{booktabs}%
\usepackage{algorithm}%
\usepackage{algorithmicx}%
\usepackage{algpseudocode}%
\usepackage{listings}
\usepackage{mathtools}   
\usepackage{physics}
\usepackage{stmaryrd}    
\usepackage{tikz-cd}     
\usepackage{esint}       
\usepackage{subcaption}  
\usepackage{orcidlink}

\usepackage[most]{tcolorbox}
\usepackage{aliascnt}
\usepackage{dsfont} 
\usepackage{tabularray}

\usepackage{geometry}
\definecolor{burgundy}{rgb}{0.5, 0.0, 0.13}
\definecolor{cinereous}{rgb}{0.6, 0.51, 0.48}
\definecolor{lightsalmonpink}{rgb}{1.0, 0.6, 0.6}

\newtheoremstyle{thmstyleone}
{18pt plus2pt minus1pt}
{18pt plus2pt minus1pt}
{\itshape}
{0pt}
{\bfseries}
{}
{.5em}
{}

\numberwithin{equation}{section}

\newtheorem{theorem}{Theorem}[section]

\newaliascnt{corollary}{theorem}
\newtheorem{corollary}[corollary]{Corollary}
\aliascntresetthe{corollary}

\newaliascnt{lemma}{theorem}
\newtheorem{lemma}[lemma]{Lemma}
\aliascntresetthe{lemma}

\theoremstyle{definition}
\newaliascnt{definition}{theorem}
\newtheorem{definition}[definition]{Definition}
\aliascntresetthe{definition}

\newaliascnt{proposition}{theorem}
\newtheorem{proposition}[proposition]{Proposition}
\aliascntresetthe{proposition}

\newaliascnt{remark}{theorem}
\newtheorem{remark}[remark]{Remark}
\aliascntresetthe{remark}

\newaliascnt{notation}{theorem}
\newtheorem{notation}[notation]{Notation}
\aliascntresetthe{notation}

\newaliascnt{assumptions}{theorem}
\newtheorem{assumptions}[assumptions]{Assumptions}
\aliascntresetthe{assumptions}

\newaliascnt{example}{theorem}
\newtheorem{example}[example]{Example}
\aliascntresetthe{example}

\newcommand{\Eq}[2]{\begin{equation}\label{#1}\begin{aligned}#2 \end{aligned}\end{equation}}

\newcommand{\theo}[2]{\rbox{\begin{theorem}\label{#1} #2 \end{theorem}}}
\newcommand{\coro}[2]{\rbox{\begin{corollary}\label{#1} #2 \end{corollary}}}

\newcommand{\prop}[2]{\bbox{\begin{proposition}\label{#1} #2 \end{proposition}}}
\newcommand{\rem}[2]{\begin{remark}\label{#1} #2 \end{remark}}

\newcommand{\ex}[2]{\begin{example}\label{#1} #2 \end{example}}

\newcommand{\bbox}[1]{\begin{tcolorbox}[arc=0mm,oversize,colback=cinereous!3!white,colframe=cinereous!100!white]#1\end{tcolorbox}}
\newcommand{\rbox}[1]{\begin{tcolorbox}[arc=0mm,oversize,colback=purple!3!white,colframe=purple!100!white]#1\end{tcolorbox}}

\renewcommand{\geq}{\geqslant}
\renewcommand{\leq}{\leqslant}

\newcommand{\eins}{\mathds{1}}
\DeclareMathOperator{\erfc}{erfc}

\title{\textbf{Partition function of 2D Coulomb gases with radially symmetric potentials and a hard wall}}

\author*[1,2]{\fnm{Matthias} \sur{Allard}\orcidlink{0000-0002-5682-424X}}\email{m.allard@unimelb.edu.au}
\author*[1]{\fnm{Peter J.} \sur{Forrester}\orcidlink{0000-0002-9328-1387}}\email{pjforr@unimelb.edu.au}
\author*[1,2]{\fnm{Sampad} \sur{Lahiry}}\email{sampad.lahiry@kuleuven.be} 
\author*[1]{\fnm{Bojian} \sur{Shen}\orcidlink{0000-0003-2467-8775}}\email{bojian.shen@unimelb.edu.au}

\affil[1]{\orgdiv{School of Mathematics and Statistics}, \orgname{University of Melbourne}, \orgaddress{\street{813 Swanston Street}, \city{Parkville, Melbourne}, \postcode{3010}, \state{Victoria}, \country{Australia}}}
\affil[2]{\orgdiv{Department of Mathematics}, \orgname{Katholieke Universiteit Leuven}, \orgaddress{\street{Celestijnenlaan 200 B bus 2400}, \city{3001 Leuven},  \country{Belgium}}}

\begin{document}

\abstract{The large $N$ asymptotic expansion of the partition function for the normal matrix model is predicted to have special features inherited from its interpretation as a two-dimensional Coulomb gas. However for the latter, it is most natural to include a hard wall at the boundary of the droplet. We probe how this affects the asymptotic expansion in the solvable case that the potential is radially symmetric and the droplet is a disk or an annulus. We allow too for the hard wall to be strictly inside the boundary of the droplet. It is observed the term of order $\log N$, has then a different rational number prefactor to that when the hard wall is at the droplet boundary. Also found are certain universal (potential independent) numerical constants given by definite integrals, both at order $\sqrt{N}$, and in the constant term. }



\maketitle
\tableofcontents

\section{Introduction}
\subsection{Confined plasmas and constrained droplets}\label{Confined plasmas and constrained droplets}
At the beginning of the 1980's Alastuey and Jancovici \cite{Alastuey1981} introduced an exactly solvable two-dimensional one-component plasma model on a disk, albeit at just one value of the dimensionless inverse temperature $\beta = 2$; this is detailed  too in the text \cite[\S 15.3.1]{Forrester2010} and the review \cite{Forrester1998}. The model is defined within a disk of radius $R$. Inside the disk is a smeared out uniform background of charge density $\rho_{\rm b}:=-N/\pi R^2$ and $N$ positive unit charged particles --- in particular the total charge is zero. The pair potential used to compute the energies associated with background and the particles is $\Psi(\vec r, \vec r') = - \log |\vec{r} - \vec{r}\,'|$ as required for a Coulomb system in two-dimensions. For general dimensionless inverse temperature $\beta$, the Boltzmann factor for this system, defined in terms of the energies $U$ of the particle-particle, particle-background and background-background interactions is (see e.g.~\cite[Eq.~(1.72) with $\Gamma = \beta$]{Forrester2010})
\Eq{0.1}{
e^{-\beta U} := e^{-\beta N^2 ((1/2) \log(R)- 3/8)} e^{- \pi \beta \rho_{\rm b} \sum_{j=1}^N | z_j|^2/2}
\prod_{1 \le j < k \le N} | z_k - z_j|^\beta,
}
with $z_j = x_j + i y_j$, $(x_j,y_j)$ being the coordinates in the plane of particle $j$.
In the case where $\rho_{\rm b}$ is fixed and
$\beta = 2$, it was shown that  the thermodynamic limit $N,R \to \infty$ of the dimensionless free energy
\Eq{0.2}{
\beta F |_{\beta = 2} := 
-\log \frac{1 }{ N!} \int_{\Omega} d^2 z_1 \cdots \int_{\Omega} d^2 z_N e^{-\beta U} \Big |_{\beta = 2} =
N^2 \log(R)- \frac{3 }{ 4} N^2 - \log \frac{\mathcal Q_N }{ N!},
}
can be obtained exactly. Here, $\Omega = \{ z: |z| <R \}$ (i.e.~the disk of radius $R$), and $\mathcal Q_N$ denotes the configuration integral
\Eq{0.2a}{
\mathcal Q_N := \int_{\Omega} d^2 z_1 \cdots \int_{\Omega} d^2z_N \, 
e^{- \pi  \rho_{\rm b} \sum_{j=1}^N | z_j|^2}
\prod_{1 \le j < k \le N} | z_k - z_j|^2.
}
Thus, the results of \cite{Alastuey1981} (at leading order in $N$), supplemented by results from \cite{Jancovici1994} (next two terms)
give the asymptotic expansion
\Eq{0.3}{
\beta F |_{\beta = 2} \sim \pi R^2 f_\Omega + (2 \pi R) \gamma_{\partial \Omega} + \frac{1 }{ 12} \log(N)+ {\rm O}(1),
}
where $f_\Omega$ (the free energy per volume) and $\gamma_{\partial \Omega}$ (the surface tension per boundary length) are given by
\Eq{0.4}{
f_\Omega = \frac{\rho_{\rm b} }{ 2} \log (\frac{\rho_{\rm b}}{2 \pi^2}), \qquad
\gamma_{\partial \Omega} = - \Big ( \frac{\rho_{\rm b} }{ 2 \pi} \Big )^{1/2} \int_0^\infty \log
\Big ( \frac{1 }{ 2} \left(1 + {\rm erf} (y)\right) \Big ) \, dy.
}

The term $\frac{1 }{ 12} \log(N)$ in (\ref{0.3}) is of much significance in the theoretical understanding of the two-dimensional one-component plasma. In fact for general $\beta > 0$ and general two-dimensional geometries, the prediction of \cite{Jancovici1994}---based on a field theory viewpoint of the electric potential and electric field---is that, for large $N$,
\Eq{0.5}{
\beta F \sim A N + B N^{1/2} + \frac{\chi }{ 12} \log(N)+ \cdots
}
where $A = \frac{\beta }{ 4 } \log (\pi \rho_{\rm b}) + f(\beta)$ for some $f(\beta)$, $B$ is proportional to the surface tension   $\gamma_{\partial \Omega}(\beta)$,
and $\chi$ the Euler characteristic of the two-dimensional domain
 (e.g.~$\chi = 1$ for a disk, $\chi = 2$ for a sphere, $\chi=0$ for an annulus). The striking feature of (\ref{0.5}) is the predicted universality of the $\log(N)$ term, being independent of $\beta$ but rather distinguishing the underlying topology. The expansion (\ref{0.5}) is consistent with the exact result (\ref{0.4}), as well as with other known exact results for the sphere \cite{Caillol1981} and annulus \cite{Fischmann2011,Fischmann2013} domains in the case $\beta = 2$. It also agrees with numerical results for the disk and sphere domains for $\beta = 4, 6$, and $8$ \cite{Tellez1999,Tellez2012}. In the case of the disk and annulus, it is expected that the structure of the expansion holds true upon the removal of the confining hard wall implied by $\Omega$, with the background uniform charge density creating an effective disk or annulus --- see the text in the following paragraph. This scenario was, in fact, the one tested numerically for the disk in \cite{Tellez1999,Tellez2012}.

 A feature of the Boltzmann factor \eqref{0.1} is that, up to a renormalisation of $\rho_{\rm b}$, its dependence on the particle coordinates is unchanged by the scalings $z_l \mapsto c z_l$
 ($l=1,\dots,N$) for any $c > 0$. Making this change of variable (with $c = \sqrt{N/\pi \rho_{\rm b}}$)  allows \eqref{0.1}, in the case $\beta = 2$, to be identified as an example of the functional form
 \Eq{2.0}{
e^{-N \sum_{j=1}^N Q(z_j)} 
 \prod_{1 \le j < k \le N} | z_k - z_j|^2,
}
specifically with $Q(z) = |z|^2$.
In fact, in this case and up to a normalisation, (\ref{2.0}) corresponds to the eigenvalue probability density function for $N \times N$ standard complex Gaussian matrices (also referred to as matrices from the Ginibre unitary ensemble, GinUE
\cite{Byun2024a}), scaled by the factor $1/\sqrt{N}$.
For more general potentials $Q$, the density \eqref{2.0} can be normalised and interpreted as the eigenvalue probability density function for the class of complex normal random  matrices $M M^\dagger = M^\dagger M$ distributed according to the probability density 
$e^{- N {\rm Tr} \, Q(M)}/Z_N$, where $Z_N$ is the corresponding normalisation constant; see the review
\cite{Zabrodin2006}. Note, however, that while the eigenvalues $\{z_j \}_{j=1}^N$ of the aforementioned normal matrix model are supported on the whole complex plane $\mathbb C$, in the model (\ref{0.1}), they are constrained to a finite disk by a hard wall boundary condition. 

Another remarkable feature of the density (\ref{2.0}) with the choice $Q(z) = |z|^2$ is that the marginal probability density on one eigenvalue (referred to as $1$-point function or level density) converges, as $N\to \infty$, to the uniform distribution $\frac{1}{\pi} \eins_{|z| < 1}$. This is known in the random matrix literature as the circular law and the support of the limiting probability measure---here the unit disk---is referred to as the droplet; see the review \cite{Bordenave2012}. 
In the rescaled variables, the (abrupt) boundary of the support $|z|=1$ is precisely the boundary of the confining disk relating to the constrained model (\ref{0.1}).

Focusing on (\ref{2.0}) with each $z_l \in \mathbb C$, one can raise it to the power $\beta/2$ to obtain a generalisation for any inverse temperature $\beta>0$. To render it a probability density, one needs to introduce the corresponding normalisation constant---partition function---${Z}_{N,\beta}^{\rm s}$ (here  the superscript ``s'' stands for soft indicating the absence of a hard wall container, and thus implying that particles are free to occupy regions outside of the droplet, albeit with low probability), explicitly given by
 \Eq{2.1}{
 {Z}_{N,\beta}^{\rm s} := \frac{1 }{ N!}
 \int_{\mathbb C} d^2 z_1 \cdots \int_{\mathbb C} d^2 z_N \,
 e^{-\beta N \sum_{j=1}^N Q(z_j)/2} \prod_{1 \le j < k \le N} |z_k - z_j|^\beta.}
It is a relatively recent result \cite{Leble2017}
(see too \cite{Sandier2015}) that this partition function, has the large $N$ asymptotic expansion
\Eq{2.1x}{
- \log {Z}_{N,\beta}^{\rm s} \sim \frac{\beta }{ 2} N^2  I_\Omega[\mu_Q] + \Big (1- \frac{\beta }{ 4}
  \Big ) N \log(N)+
 \bigg ( C(\beta) + \Big (1- \frac{\beta }{ 4}
  \Big )E_\Omega[\mu_Q] \bigg ) N + {\rm o}(N).
 }
 Here $I_\Omega[\mu]$ and $E_\Omega[\mu]$ are, respectively, the energy and entropy functionals specified by
 \Eq{I1}{I_\Omega[\mu]= -\int_{\Omega^2} \log | z - w| \, d\mu(z) d\mu(w)+ \int_{\Omega}Q(z) \, d\mu(z), \quad E_\Omega[\mu]= \int_{\Omega} \mu(z) \log( \mu(z))\, d^2\mu(z),
 }
 respectively.
The measure $\mu_Q$ --- referred to as the equilibrium measure and supported on $\Omega$ --- is chosen to minimise $I_\Omega[\mu]$, while $C(\beta)$ is not known in general but is independent of $Q$. In the case of that the equilibrium measure is flat, the best known bound on the remainder in (\ref{2.1x}) is O$(N^{1/2} \log N)$ \cite{Armstrong2021}, and thus close to the predicted O$(\sqrt{N})$ in  keeping with (\ref{0.5}). However without this assumption, rigorous bounds on the remainder are of higher orders \cite{Bauerschmidt2019}, \cite{Serfaty2023}.

We see that the scaled configuration integral $\mathcal Q_N/N!$ as specified by (\ref{0.2a}) with $\rho_b = N/\pi$ is the case $\beta = 2$,
$Q(z) = |z|^2$ of (\ref{2.1}), yet with $\mathbb{C}$ replaced by $\Omega$. Nonetheless, making use of (\ref{0.2}), (\ref{0.3}) with $R=1$, and (\ref{0.4}) 
 we can check that the large $N$
expansion up to order $N$ agrees with (\ref{2.1}) (this requires the computation that 
$I_\Omega[\mu_Q] = -3/4$ for $\mu_Q$ the circular law). However, this agreement only partially persists at higher order. Specifically, the next two terms in (\ref{2.1x}) with $Q(z) = |z|^2$ and general $\beta > 0$ are conjectured to be
\Eq{2.2}{a_\beta \sqrt{N} + \frac{1 }{ 12}\log(N),\qquad a_\beta = \frac{4 \log(\beta/2) }{ 3 \pi^{1/2}}.}
Here the value of $a_\beta$ is a conjecture
reported in
\cite[Eq.~(3.2)]{Can2015}, \cite[Eq.~(4.7)]{Byun2024a} (note that this vanishes for $\beta = 2$), and the reasoning for the $\log(N)$ term the same as for (\ref{0.5}) in the case of a disk. For 
$\beta = 2$ the conjectured values can be checked using the exact result 
$ {Z}_{N,\beta}^{\rm s} |_{\beta = 2}
= \pi^N N^{-N(N+1)/2} \prod_{j=1}^{N-1} j!$ (equivalent to (\ref{1.18}) below with $q(r) = r^2$ and the upper terminal of integration now $\infty$). To contrast with  (\ref{0.3}), note that the $\log(N)$ term agrees, while there is now a non-zero coefficient of $\sqrt{N}$.
 
The similarities and differences between the partition function asymptotics of the droplet model in the case $Q(z) = |z|^2$ without and with a hard wall boundary
suggests an investigation of this question for more general potentials.
To maintain  the property that the partition function (\ref{2.1}) can be computed exactly we require $\beta =2$,  
$Q(z) = q(r)$, $r=|z|$ (radial symmetry of the potential, which will further be taken to be globally sub-harmonic, and strictly subharmonic on the droplet),     and that the hard wall container be a disk about the origin. 
Let us denote the partition function (\ref{2.1}) in this setting by $Z^{\rm s}_N[q(r)]$, and its companion with the hard wall by  $Z^{\rm h}_N[q(r)]$.

Here the hard wall boundary defined by the confining unit disk need not correspond to a boundary of the droplet (from our assumption on the form of the potential, the droplet may be either a disk or an annulus). This circumstance  is well motivated from the consideration of the probabilities $\{E_{N}(n;|z|<R)\}$ that there are exactly $n$ of the $N$ particles inside the disk of radius $R$ for the Coulomb gas system corresponding to (\ref{2.1}) with $\beta = 2$; see e.g.~\cite[\S 3.1]{Byun2024a} and \cite{Charlier2023},
\cite{Byun2024} for recent developments with $n=0$.  Relevant to the present study is the particular choice
$n=N$ when one has
\Eq{2.2b}{
E_{N}(N;|z|<1) = \frac{{Z}_{N}^{\rm h}[q(r)] }{ {Z}_{N}^{\rm s}[q(r)]}.}
In the case $q(r) = (s r)^2$, ($0< s < 1$) the asymptotic expansion of $\log E_{N}(N;|z|<1)$ is known
(\cite[first 3 terms]{Cunden2016}, \cite[Th.~1.4 with $g=1,\alpha=0,b=1$]{Charlier2024}) as is the asymptotic expansion of
$\log {Z}_{N}^{\rm s}$ (recall the text two sentences below (\ref{2.2})). Hence
\Eq{EN}{\log \bigg ( \frac{{Z}_N^{\rm h}((s r)^2) }{ (2 \pi)^N} \bigg )  = \tilde{C}_1 N^2 + \tilde{C}_2 N \log(N)+ \tilde{C}_3N + \tilde{C}_4 \sqrt{N} - \frac{1 }{ 3} \log(N)+ {\rm O}(1) ,}
where
\Eq{}{\tilde{C}_1 = \frac{s^4 }{ 4} - s^2   , \quad \tilde{C}_2 = \frac{s^2 -2 }{ 2} , \quad \tilde{C}_3 = (1 - s^2) \Big ( 1 - \log  (1 - s^2) \Big ) - s^2 \log (s/\sqrt{2\pi}) ,\nonumber}
\Eq{EN1}{
\tilde{C}_4 = \sqrt{2} s \bigg ( \int_{-\infty}^0 \log \Big ( \frac{1 }{ 2} \erfc(y) \Big ) \, dy +
\int_0^\infty  \log \Big ( \sqrt{\pi} y e^{y^2} \erfc(y) \Big )  \, dy \bigg ),}
with the O$(1)$ constant term also available in \cite{Charlier2024}.

In the setting of (\ref{EN}),
the equilibrium measure for the hard wall system consists of a uniform bulk density $\frac{s^2 }{ \pi }$ for $|z|<1$, and a uniform surface density of total mass $(1-s^2)$ on the boundary $|z|=1$ --- it is an example of a balayage measure \cite{Ameur2024a}. From this, the energy  functional in (\ref{I1}) can be evaluated and substituted in (\ref{2.1x}) to reclaim $\tilde{C}_1$. Also,
the first integral in the formula for $\tilde{C}_4$ can be rewritten as the integral in (\ref{0.4}), and one notes the appearance of a simple fraction as the  coefficient of the $\log(N)$ terms, although different from that in
(\ref{2.2}). In the present study we can probe the universality of this coefficient with respect to the topology of the droplet and the external potential. We draw attention to the recent work \cite{Ameur2024} for a study of correlations near the hard wall for analogous hard wall problems.

In addition to the topological significance of the coefficient of the $\log(N)$ term in (\ref{2.2}) as it applies to the asymptotic expansion of $\log  \tilde{Z}_{N}^{\rm h}$, since the work of Wiegmann and Zabrodin \cite{Zabrodin2006a} (see also \cite{Klevtsov2019}) there is known to be topological/geometrical significance in the constant term of the large $N$ expansion of the latter, conjectured to contain an additive term that can be written in terms of the  regularised determinant of the Laplacian; see \cite[\S 9.3]{Serfaty2024} for a clear statement.
We remark that
 a check beyond  that of a rotationally symmetric droplet has recently been found in \cite{Byun2025c}.
The structure of the constant term has been further developed in \cite{Byun2023}, where it is predicted 
 that for a  droplet $\mathbb S$ there will be a further  additive term at O$(1)$ involving the derivative of the Riemann zeta function and the Euler characteristic, 
\Eq{2.2c}{ \chi \zeta'(-1);} see also \cite[Eq.~(5.17)]{Byun2024a}.
Since our asymptotic expansions for
$\log {Z}_{N}^{\rm h}[q(|z|)]$ extend to the constant term, we have the opportunity to further test this.

\subsection{Description of the model}
\begin{figure}[H]
        \centering
        \includegraphics[width=1\textwidth]{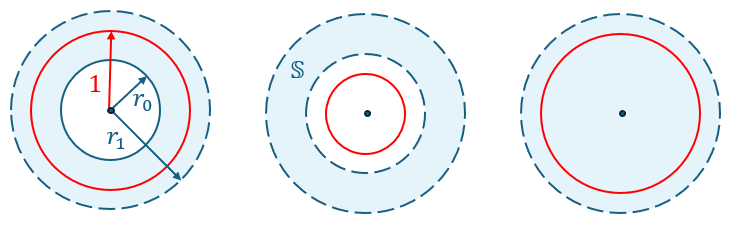}
        \captionsetup{justification=justified}
        \caption{\label{fig: config} Hard wall and droplet configurations: \emph{(Left)} $0<r_0<1\leq r_1$: The droplet (in shaded blue) is an annulus. The hard wall (red) is placed inside the droplet or exactly at $r_1$.
        \emph{(Middle)} $1\leq r_0< r_1$: The droplet is an annulus. The hard wall is placed inside the inner radius of the droplet or exactly at $r_0$. \emph{(Right)} $1\leq  r_1$: The droplet is a disk. The hard wall is placed inside the droplet or exactly at $r_1$. Note: The boundaries drawn here are distinct from  those corresponding to the equilibrium measure of the confined Coulomb gas, which  have outer boundary $r=1$, inner boundary $r_0,1,0$ respectively.}   
\end{figure}

As previously said, we are considering the 
$\beta = 2$ case of the partition function (\ref{2.1}) with each integration domain the unit disk, thought of as a confining domain. The one-body potential is taken to be radially symmetric $Q(z) = q(r)$, $r = |z|$, with the further requirement that it is strictly subharmonic,
\Eq{eq: Laplacian}{
\Delta q(r)=\frac{1}{r}(rq'(r))'=q''(r)+\frac{q'(r)}{r}>0.
}
Standard theory \cite{Saff1997}, conveniently summarised in \cite[\S 1.2]{Byun2023}, gives that upon the assumption of lower bound on the growth of $q(r)$, without confining domain,
the equilibrium measure \(\mu_Q\) associated with the potential \(Q\) is specified by
\Eq{eq: droplet}{
d\mu_Q = \frac{1}{4}\Delta Q \cdot \eins_{\mathbb{S}}\, dA,\quad dA(z):=\frac{d^2z}{\pi}, \quad
\mathbb{S} := \mathbb{A}_{r_0, r_1} = \{ z \in \mathbb{C} \mid r_0 \leq |z| \leq r_1 \},
}
where $\eins_{\mathbb{S}}$ is the indicator function of the set $\mathbb{S}$---the droplet.
Here \(r_0\), the inner radius of the droplet, is the smallest non-negative solution to the equation \(q'(r_0) = 0\). In particular, the droplet becomes a disk (i.e., \(r_0 = 0\)) if \(q'(0) = 0\), and an annulus otherwise. The outer radius \(r_1\) is the unique solution to
$
r_1 q'(r_1) = 2$.
Superimposing now the hard wall of the confining domain  $\mathbb{D}:=\{z\in\mathbb{C}\mid \abs{z}\leq1\}$ leads to several interesting configurations, depending on the values of \(r_0\) and \(r_1\) which are illustrated in \autoref{fig: config}.

\rem{}{
Due to the presence of the hard wall, in the middle setting of 
\autoref{fig: config} it is not
necessary that the droplet has finite boundaries, i.e.~$r_0, r_1$ may be $\infty$. It is then valid to consider the non-confining sub-harmonic potential $q(r)=  c / r$, $c > 0$, for example. 
(If we were to drop the requirement that $q(r)$ being sub-harmonic, taking $q(r) = 0$ similarly corresponds to $r_0, r_1$ being infinite; we make mention of this case due to the recent work on the asymptotics of the corresponding partition function when the hard wall boundary is a general Jordan curve with corners \cite{Johansson2023}.)
For the  left and middle settings of \autoref{fig: config}, it is also permitted that $r_1$ equal $\infty$. An example is for the potential $-c \log( r) +\log(1+r^2)$ as is encountered in studies of the Coulomb gas on the sphere, stereographically projected to the plane \cite{Fischmann2011,Byun2025b,Byun2025a,Byun2025}.}

Consider now the partition function ${Z}_{N}^{\rm h}[q(r)]$
as defined at the end of the paragraph above (\ref{2.2b}).
A standard calculation \cite[Exercises 15.3 q.1(ii)]{Forrester2010} gives that the corresponding $N$-dimensional integral factorises as the product of $N$ one-dimensional integrals. Specifically, after taking the logarithm,
\Eq{1.18}{
\log(\frac{{Z}_N^{\rm h}[q(r)]}{(2 \pi)^N })= \sum_{j=0}^{N-1}\log u_{j}, \quad u_j:=\int_{0}^{1}r^{2j+1}e^{-N q(r)}dr.
}
The essential mathematical question then is to have sufficient control of the large $N$ asymptotic of the integral, to be carried out using Laplace's method, so that the asymptotics of the sum can then be determined using the Euler-Maclaurin summation formula. In the case equivalent to choosing $q(r) = \alpha r^2$, $0<\alpha<1$, 
and requiring that all the particles lie {\it outside} of a hard wall boundary at $r=1$ (this shows itself in the computation of the gap probability of no eigenvalues in $\mathbb D_{\sqrt{\alpha}}$ for scaled GinUE), such a strategy was first implemented in \cite{Forrester1992} to obtain an expansion of the form of (\ref{EN}), with the first four coefficients determined explicitly.
However in that case the integral corresponding to $u_j$ in (\ref{1.18}) --- which now has terminals of integration $[1,\infty)$ --- can be expressed in terms of the incomplete gamma function. There is no need to apply Laplace's method to determine the asymptotics as this is available in the literature (see e.g.~\cite[\S 8.11]{DLMF}). This is similarly true of the choice of potential $q(r) = r^{2b} + (2a/N) \log(r)$ considered in the revival due to Charlier \cite{Charlier2024} of the study of asymptotic expansions relating to the study \cite{Forrester1992} (due to the dependence on $N$ of the coefficient $\log(r)$, this choice falls outside of the class permitted in the present work, except for $a=0$. The recent study of Byun et al.~\cite{Byun2023} on the asymptotic expansion of (\ref{2.1}) for $\beta = 2$ and a radially symmetric potential $Q(z) = q(|z|)$, which thus is without a confining hard wall boundary (it is this study which has directly motivated the present work) has both the complexities of requiring first principal Laplace type asymptotic expansions, and the use of the Euler-Maclaurin summation. See too the subsequent works \cite{Ameur2023} (allowing for the droplet to be disconnected) and
\cite{Byun2025b} ($\beta = 2$ Coulomb gas on a sphere).

\section{Main results}
In the left and right cases of \autoref{fig: config} the effect of the hard wall on the equilibrium measure is to restrict $\mu_Q$ as specified in (\ref{eq: droplet}) to $\mathbb{A}_{r_0,1}=\mathbb{S}\cap \mathbb{D}$, and furthermore to induce a uniform surface charge of total mass
\Eq{eq: def eta}{
\eta:=1- \int_{A_{r_0,1}} d \mu_Q = 1 - \frac{1 }{ 2} \int_{r_0}^1 (r q'(r))' \, dr = 1 - \frac{ q'(1) }{ 2};
}
see e.g.~\cite{Ameur2024a}.
If $\eta=1$, the hard wall is exactly at $r_0$. For $\eta\in]0,1[$, the hard wall is between $r_0$ and $r_1$. If $\eta=0$, the hard wall is exactly at $r_1$. Note that this quantity is still well defined in the situation of the middle figure, when it satisfies
$\eta>1$. In this circumstance, the equilibrium measure in the presence of the hard wall is uniform on the circle $|z|=1$. 

We introduce some universal (i.e.~$q(r)$ independent) integrals appearing in the expansions, which in their integrands typically involve
\Eq{eq: def phi}{
\Phi(x):=\log\left(\frac{1}{2}\erfc(x)\right),
\quad \Phi'(x)=-\frac{2 e^{-x^2}}{\sqrt{\pi } \erfc(x)
},
}
where $\erfc$ is the usual complementary error function.
Thus we set
\Eq{eq: pos constant in}{
\alpha^{\rm in}&:=\int_{-\infty}^0  \frac{ \left(x^2-2\right)}{3 }\Phi'(x) dx\approx 0.36941\\
\beta^{\rm in}&:=-\int_{-\infty}^0  \frac{\left(x^2+1\right) }{6}\Phi'(x) dx\approx 0.16186\\
\gamma^{\rm in}&:=-\frac{1}{\sqrt{2}}\int_{-\infty}^0 \Phi(x) dx\approx 0.23876,
}
and 
\Eq{eq: pos constant out}{
\alpha^{\rm out}&:=\int_0^{\infty}\left(  \frac{1}{3}\left[x(2x^2-3)+(x^2-2)\Phi'(x) \right]-\frac{1}{x + 1}\right)dx\approx 0.27752\\
\beta^{\rm out}&:=-\frac{1}{6}\int_0^{\infty}\left[x(2x^2+3)+(x^2+1)\Phi'(x) \right]dx\approx 0.14742\\
\gamma^{\rm out}&:=-\frac{1}{\sqrt{2}}\int_0^{\infty} \log\left( \sqrt{\pi} x e^{x^2} \erfc(x) \right) dx\approx 0.91194.
}
In (\ref{2.23}) below we give equations implying that $\beta^{\rm in}$ ($\beta^{\rm out}$) can be written in terms of $\alpha^{\rm in}$ ($\alpha^{\rm out}$).

Also required are
the logarithmic energy and entropy of $\mu_Q$ 
(\ref{eq: droplet}) restricted to
$\mathbb{S}\cap \mathbb{D}=\mathbb{A}_{r_0,1}$,
\Eq{eq: energ muQ}{
I_{\mathbb{S}\cap \mathbb{D}}[\mu_{Q}]:=&\,\frac{1}{2}\int_{\mathbb{S}\cap \mathbb{D}}Q(z)d\mu_Q(z)+\frac{1}{2}q(1)(1-\eta)
=\,\frac{1}{8}\int_{\mathbb{S}\cap \mathbb{D}}Q(z)\Delta Q(z)dA(z)+\frac{1}{4}q(1)q'(1)\\
 \quad =&\,\frac{1}{2}\int_{r_0}^1u \Delta Q(u) \big[q(u)-uq'(u)\log(u) \big]  du,\\
}
\Eq{eq: ent muQ}{
E_{\mathbb{S}\cap \mathbb{D}}[\mu_{Q}]:=\,\int_{\mathbb{S}\cap \mathbb{D}}\log(\frac{1}{4}\Delta Q(z))d\mu_{Q}(z)
=\,\frac{1}{2}\int_{r_0}^1\log(\frac{1}{4}\Delta Q(u))u\Delta Q(u)du,}
as well as
\Eq{eq: F annulus}{
 F_{\mathbb{S}\cap \mathbb{D}} [Q]:=\frac{1}{12}\log (\frac{r_0^2\Delta Q(r_0)}{\Delta Q(1)})-\frac{1}{16}\left(\frac{\partial_r \Delta Q(1)}{\Delta Q(1)}-r_0\frac{\partial_r \Delta Q(r_0)}{\Delta Q(r_0)}\right)+\frac{1}{24}\int_{r_0}^{1}\left(\frac{\partial_r \Delta Q(r)}{\Delta Q(r)}\right)^2rdr
}
\Eq{eq: F disk}{
 F_{\mathbb D} [Q]:=\frac{1}{12}\log (\frac{4}{\Delta Q(1)})-\frac{1}{16}\frac{\partial_r \Delta Q(1)}{\Delta Q(1)}+\frac{1}{24}\int_{0}^{1}\left(\frac{\partial_r \Delta Q(r)}{\Delta Q(r)}\right)^2rdr.
}
We remark that a generalisation of the latter two quantities first appeared in \cite[Eqns.~(1.16) and (1.22) with slightly different notation]{Byun2023} involving an extra parameter $r_1$ --- setting $r_1 = 1$ gives (\ref{eq: F annulus}) and (\ref{eq: F disk}). Note too that the assumption that the potential is strictly subharmonic on the droplet, $\Delta Q(r) > 0$, guarantees that the integrals in these functionals are well defined.

In terms of these integrals and functionals, we can now state our findings for the asymptotics of the logarithm of the partition function ${Z}_N^{\rm h}[q(r)]$
as given by the sum (\ref{1.18}).

\subsection{In/Out Annulus case: $0 < r_0 < 1\leq r_1$}
\theo{theo: In/Out Annulus}{
\textup{(In/Out Annulus case)} 
Consider the setting of the left diagram in
\autoref{fig: config}. 
The large $N$ expansion of the logarithm of the partition function ${Z}_N^{\rm h}[q(r)]$ depends on $r_1$ being strictly greater than 1, or $r_1=1$. Specifically, with $\mathbb S$ denoting the droplet and $\mathbb D$ denoting the unit disk,
\begin{itemize}
 \item If $r_1>1$ $(1>\eta>0)$: 
\Eq{}{
& \log(\frac{{Z}_{N}^{\rm h}[q(r)]}{ (2\pi)^N })= \,- N^2 \left[I_{\mathbb{S}\cap \mathbb{D}}[\mu_{Q}]+\eta q(1) \right]- \frac{(1+\eta)}{2} N\log(N)\\
& \quad - N\left[\frac{1}{2} E_{\mathbb{S}\cap \mathbb{D}}[\mu_{Q}]  - (1-\eta) \frac{\log(\pi/2)}{2}  + \eta (\log (2 \eta)-1) \right]
- \sqrt{N}(\gamma^{\rm in}+\gamma^{\rm out})\sqrt{\Delta Q(1)}- \frac{\log(N)}{4}\\
& \quad +F_{\mathbb{S}\cap \mathbb{D}}[Q] -(\alpha^{\rm in}+\alpha^{\rm out})+(\beta^{\rm in}+\beta^{\rm out})  \frac{\partial_r \Delta Q(1)}{\Delta Q(1)} 
+\frac{\Delta Q(1)}{4 \eta }-\frac{1}{2} \log(\eta) +\frac{1}{4} \log(2\pi  \Delta Q(1))+{\rm o}(1). 
}
\item If $r_1=1$ $(\eta=0)$: 
\Eq{}{
&\log(\frac{{Z}_{N}^{\rm h}[q(r)]}{ (2\pi)^N})= \,- N^2 I_{\mathbb{S}\cap \mathbb{D}}[\mu_{Q}]- \frac{1}{2}N\log(N)- N\left[\frac{1}{2} E_{\mathbb{S}\cap \mathbb{D}}[\mu_{Q}]  - \frac{\log(\pi/2)}{2}  \right]\\
& \qquad - \sqrt{N}\gamma^{\rm in}\sqrt{\Delta Q(1)}+F_{\mathbb{S}\cap \mathbb{D}} \left[ Q \right]-\alpha^{\rm in}+\beta^{\rm in}  \frac{\partial_r \Delta Q(1)}{\Delta Q(1)}+ \frac{\log(2)}{2}+{\rm o}(1). 
}
\end{itemize}
The positive universal constants  $\alpha^{\rm out/in}, \beta^{\rm out/in}, \gamma^{\rm out/in}$ are given by \eqref{eq: pos constant in} and \eqref{eq: pos constant out}. The functionals $I_{\mathbb{S}\cap \mathbb{D}}[\mu_{Q}]$,  $E_{\mathbb{S}\cap \mathbb{D}}[\mu_{Q}]$, $F_{\mathbb{S}\cap \mathbb{D}}\left[ Q \right]$ are respectively given in \eqref{eq: energ muQ}, \eqref{eq: ent muQ} and \eqref{eq: F annulus}.
}

\ex{E1}{\normalfont ${}$ \\
1.~Consider the particular radial potential
\Eq{q2}{q(r) = r^2 - 2 a \log(r), \qquad 0<a<1.
}
From the theory below (\ref{eq: droplet}) 
\Eq{q2a}{r_0 = \sqrt{a}, \quad r_1 = \sqrt{a+1}, \quad \eta = a.}
Hence for $a$ in the interval $]0,1[$ as specified, the setting of the first subcase of
\autoref{theo: In/Out Annulus} is obtained.
After substituting (\ref{q2}) in the definition (\ref{1.18}) of $u_j$, a simple change of variables gives the special function evaluation
\Eq{q3}{
u_j = \frac{1 }{ 2} N^{-(j+aN + 1)} \gamma(j+aN+1;N),
}
where $\gamma(a;z) := \int_0^z t^{a-1} e^{-t} \, dt$ denotes the lower incomplete gamma function.\\
2.~To obtain an example of the second subcase of \autoref{theo: In/Out Annulus}, consider the potential
\Eq{q4}{q(r) = c r^2 - 2 (c - 1) \log(r), \qquad c>1.}
For this
\Eq{q4a}{r_0 = \sqrt\frac{c-1 }{ c}, \quad r_1 = 1, \quad \eta = 0.}
Furthermore
\Eq{q4b}{
u_j = \frac{1 }{ 2} (Nc)^{-(j+(c-1)N + 1)}
\gamma(j+(c-1)N+1;cN).
}
3.~As the incomplete gamma function can be computed to arbitrary position in computer algebra software such as Mathematica, we have carried out numerical checks of \autoref{theo: In/Out Annulus} in these cases, with agreement found.
}

\subsection{In/Out Disk case: $0 = r_0 < 1\leq r_1$}
\theo{theo: In/Out Disk}{
\textup{(In/Out Disk case)} 
Consider the setting of the right diagram in
\autoref{fig: config} (denoted disk). 
The large $N$ expansion of the logarithm of the partition function ${Z}_N^{\rm h}[q(r)] |_{\rm disk}$, which depends on $r_1$ being strictly greater than 1, or $r_1=1$, is simply related to the corresponding results for $\tilde{Z}_N[q(r)] |_{\rm annulus}$ of 
\autoref{theo: In/Out Annulus} by the relation
\Eq{}{
\log(\frac{{Z}_{N}^{\rm h}[q(r)]|_{\rm disk}}{ (2 \pi)^N })=  \log(\frac{{Z}_{N}^{\rm h}[q(r)]|_{\rm annulus}}{ (2 \pi)^N }) - \frac{1 }{ 12} \log(N)+ \zeta'(-1) + F_{\mathbb D}[Q] - F_{\mathbb{S}\cap \mathbb{D}}[Q],
}
where $\zeta$ is the Riemann zeta function.
}

\ex{E2.3}{
Consider the radial potential
\Eq{}{
q(r)=a^2 r^2, \quad 0<a<1.
}
Then we have
\Eq{}{
r_0=0,\quad r_1=\frac{1}{a},\quad \eta=1-a^2,
}
which fits in the setting of \autoref{theo: In/Out Disk}. In this case, $u_j$ is given by
\Eq{}{
u_j=\frac{1}{2}(a^2N)^{-j-1}\gamma(j+1;a^2N).
}
The asymptotics in this case is the expansion (\ref{EN}), and thus known independently of our work. We have confirmed agreement with \autoref{theo: In/Out Disk} (see \S 2.4 for a discussion relating to comparing the terms at ${\rm O}(1)$ given by definite integrals).
}

\subsection{Out Annulus case: $1 \leq r_0 <r_1$}
\theo{theo: Out Annulus}{
\textup{(Out Annulus case)} 
Consider the setting of the middle diagram in
\autoref{fig: config}, so that $1 \leq r_0 <r_1$.
The large $N$ expansion of the logarithm of the partition function ${Z}_N^{\rm h}[q(r)]$ depends on $r_0$ being strictly greater than 1, or $r_0=1$. Specifically
\begin{itemize}
 \item If $r_0>1$ $(\eta>1)$: 
\Eq{2.21}{
\log(\frac{{Z}_{N}^{\rm h}[q(r)]}{ (2 \pi)^N })=&\, -N^2 q(1)-N\log(N)+N (1-\log(2)+ (\eta-1)\log(\eta-1)-\eta\log\eta)\\
&+\frac{1}{2}\log \frac{\eta-1}{\eta}+\frac{\Delta Q(1)}{4\eta(1-\eta)}+{\rm o}(1).
}
\item If $r_0=1$ $(\eta=1)$: 
\Eq{}{
\begin{aligned}
    \ 
\log(\frac{{Z}_{N}^{\rm h}[q(r)]}{ (2 \pi)^N })=&\,-N^2 q(1) - N \log(N)- N (\log(2) - 1) - \sqrt{N} \gamma^{\rm{out}} \sqrt{\Delta Q(1)} - \frac{\log(N)}{4} \\
& - \alpha^{\rm{out}} + \beta^{\rm{out}} \frac{\partial_r \Delta Q(1)}{\Delta Q(1)} + \frac{\Delta Q(1)}{4} + \frac{1}{4} \log \left( \frac{\pi \Delta Q(1)}{2} \right) + \mathrm{o}(1).
\end{aligned}
}
\end{itemize}

}

\ex{}{
The potential $q(r)$ in \autoref{E1}.1, but now with parameter chosen as $a>1$, and $a=1$, give examples for the two cases in \autoref{theo: Out Annulus} respectively.
}
\begin{remark} 
In the case where \( r_1 < 1 \), that is, when the unconstrained droplet lies entirely on the interior of the hard wall, the equilibrium configuration is unaffected by the presence of the wall. In this setting, for large $N$ one expects   
\[
\log\left(\frac{Z_{N}^{\mathrm{h}}[q(r)]}{(2\pi)^N}\right) \sim \log\left(\frac{Z_{N}^{\mathrm{s}}[q(r)]}{(2\pi)^N}\right),
\]  
where $Z_{N}^{\mathrm{s}}$ is specified by (\ref{2.1}) with $\beta = 2$. For the asymptotic expansion of this latter quantity, up to terms which go to zero with $N$, we have the known results 
 \cite[Eq.~(1.17)]{Byun2023} in the annular case, and \cite[Eq.~(1.23)]{Byun2023} in the disk case. As $r_1$ approaches the hard wall boundary from below, by analogy with terminology used previously in the case of the log-gas on the line \cite{Claeys2008}, one may say that the soft edge meets the hard edge. 
\end{remark}
\subsection{Discussion}\label{S2.4}
Highlighted in the Introduction was the universal term $(\chi/12) \log(N)$ in the free energy expansion (\ref{0.3}). It was commented that this is expected to be valid both the settings that the hard wall corresponds to the boundary of the droplet, and when the hard wall is removed altogether. Our result of \autoref{theo: In/Out Annulus}, case $r_1=1$ ($\eta = 0$), indeed confirms the former scenario for the annulus droplet, with there being no terms proportional to $\log N$. Using this fact in \autoref{theo: In/Out Disk}
then confirms that there is a term proportional to $\log N$, proportionality constant $-\frac{1 }{ 12}$ (the minus sign is due to the free energy being minus the logarithm of the partition function), again in agreement with predictions.

A striking feature of our expansions, beyond the cases covered by the conjectures from \cite{Jancovici1994}, is that the $\log(N)$ term remains universal --- being independent of the potential $q(r)$ and possibly $\beta$ (although our study is restricted to $\beta = 2$), albeit depending on whether $\eta$ is zero, is between 0 and 1, is equal to 1, or is greater than 1. Specifically, for the annulus case $r_1 > 1$ ($1 > \eta > 0$) of \autoref{theo: In/Out Annulus} there is the term
$- \frac{1 }{ 4} \log (N)$; for the disk case $r_1 > 1$ ($1 > \eta > 0$), substituting this fact in the result of \autoref{theo: In/Out Disk} shows that there is the term $-\frac{1 }{ 3} \log(N)$; and for the annulus case with $r_0>1$ ($\eta > 1$) there is no $\log(N)$ term, but for $r_0 = 1$ ($\eta = 1)$ there is the term $- \frac{1 }{ 4} \log (N)$ as for the annulus case with $r_1 > 1$ ($1 > \eta > 0$). How one might anticipate these findings remains an open question.\footnote{On this point S.~Byun (private correspondence) has drawn our attention to the work \cite{Charlier2021}. Contained therein is the asymptotic expansion of the logarithm of the partition function for the one-dimensional log-gas with Gaussian, Laguerre and Jacobi type weights (so defined by the type of singularities at the spectrum edge; i.e.~two soft edges in the Gaussian case with density decaying as a square root, one hard and one soft edge for the Laguerre case, and two hard edges for the Jacobi case). The relevant point with respect to our findings is that the coefficient of $\log(N)$ in the three cases is $-\frac{1 }{ 12}$, $-\frac{1 }{ 6}$ and $-\frac{1 }{ 4}$ respectively, and thus proportional to one plus the number of hard edges. On the other hand, in distinction to the two-dimensional case, for the log-gas, the coefficient of $\log(N)$ is not independent of $\beta$. Thus taking the Gaussian $\beta$ ensemble as an explicit example, one has the known exact result that it is given by $\frac{1 }{ 12}(\beta/2 + 2/\beta) - \frac{1 }{ 4}$ \cite[\S 4.3]{Mironov2012}.}

Focusing attention next of the term proportional to $\sqrt{N}$, present in all the large $N$ expansions except (\ref{2.21}), one encounters a proportionality constant consisting of a factor of $\sqrt{\Delta Q(1)}$, multiplied by the universal numerical constants $\gamma^{\rm out}, \gamma^{\rm in}$ or their sum. Already $\gamma^{\rm in}$ was known from the work \cite{Jancovici1994} deriving (\ref{0.3}), where it has the interpretation of a surface tension. The sum of $\gamma^{\rm out}, \gamma^{\rm in}$
was first encountered in the study \cite{Forrester1992} of the large $N$ asymptotics of the gap probability $E_N(n;|z|<\alpha \sqrt{N})$, $0 < \alpha < 1$ in the setting cited above (\ref{2.2b}). The universal constants $\alpha^{\rm in}, \alpha^{\rm out},\beta^{\rm in},\beta^{\rm out}$ at O$(1)$ are particular versions of constants first obtained in the asymptotic expansions of gap probabilities, in \cite{Charlier2024}, for Mittag-Leffler type potentials $Q(z)=\abs{z}^{2b}+\frac{2 a}{N}\log(\abs{z})$. Thus, in \cite[Th.~1.4]{Charlier2024}, the constant term comprises terms proportional to $2b=\partial_r \Delta Q(1)/\Delta Q(1)+2$ (where $b$ is the parameter of the Mittag-Leffler potential) which can then be identified as $\tilde{\beta}^{\rm in}, \tilde{\beta}^{\rm out}$. Explicitly,
\Eq{}{
& \tilde{\beta}^{\rm in} := \int_{-\infty}^0 g(y) \, dy, \qquad
g(y) := 2 y \log \Big ( \frac{1 }{ 2} \erfc(y) \Big ) + \frac{e^{-y^2} (1 - 5 y^2)}{ 3 \sqrt{\pi} \erfc(y)}, \nonumber \\
& \tilde{\beta}^{\rm out} :=\frac{1}{2}+
\int_0^\infty \bigg ( g(y) + \frac{11 }{ 3} y^3 + 2y \log y + \Big ( \frac{1 }{ 2} +
2 \log(2 \sqrt{\pi}) \Big ) y \bigg ) \, dy,
}
and it can be shown that $\tilde{\beta}^{\rm in}, \tilde{\beta}^{\rm out}$ agree with the constants $\beta^{\rm in}, \beta^{\rm out}$ in \eqref{eq: pos constant in} and \eqref{eq: pos constant out}. Indeed, using integration by parts $\int_{a}^{b}x^2 \Phi'(x)dx=x^2\Phi(x)\big|_{a}^{b}-\int_{a}^{b}2x \Phi(x)dx $, with $\Phi(x)$ defined in \eqref{eq: def phi}, we have
\Eq{}{
\tilde{\beta}^{\rm in}=\int_{-\infty}^0 \left(2x \Phi(x)+\frac{(5x^2-1)}{6}\Phi'(x)\right)dx=-\int_{-\infty}^0 \frac{(x^2+1)}{6}\Phi'(x)dx=\beta^{\rm in}.
}
On the other hand,
\Eq{}{
\tilde{\beta}^{\rm out}-\beta^{\rm out}&=\frac{1}{2}+\int_{0}^{\infty}(2x\Phi(x)+x^2\Phi'(x)+4x^3+(1+2\log(2\sqrt{\pi}))x+2x\log x)dx\\
&=\frac{1}{2}+\lim_{z\to\infty}\left(z^2\Phi(z)+\int_{0}^{z}(4x^3+(1+2\log(2\sqrt{\pi}))x+2x\log x)dx\right)=0.
}
Thus,
\Eq{2.23}{
\tilde{\beta}^{\rm in}=\beta^{\rm in}= - \frac{1}{2} \alpha^{\rm in} + \frac{\log(2)}{2} ,\qquad \tilde{\beta}^{\rm out}= \beta^{\rm out}=-\frac{1 }{ 2} \alpha^{\rm out} + \frac{1 }{ 4} \log(\pi),
}
with the second equality in each of the equations following from similar manipulations. Hence, one can use \eqref{2.23} to eliminate $\alpha^{\rm in},\alpha^{\rm in}$ in \autoref{theo: In/Out Annulus}, \autoref{theo: In/Out Disk} and \autoref{theo: Out Annulus}. In particular, we have
\Eq{}{
-(\alpha^{\rm in}+\alpha^{\rm out})+(\beta^{\rm in}+\beta^{\rm out})  \frac{\partial_r \Delta Q(1)}{\Delta Q(1)} =(\beta^{\rm in}+\beta^{\rm out}) \left( \frac{\partial_r \Delta Q(1)}{\Delta Q(1)}+2\right)-\log(2)-\frac{\log(\pi)}{2},
}
where the RHS corresponds to the formulation in  \cite[Th.~1.4]{Charlier2024}.

Still regarding the constant term, it was highlighted in \autoref{Confined plasmas and constrained droplets} that  \cite{Byun2023} predicts the universal additive term (\ref{2.2c}) distinguishing the disk and annulus geometry. This is consistent with our findings in \autoref{theo: In/Out Annulus} and \autoref{theo: In/Out Disk}, as it remains unaffected by the hard wall being strictly inside of the droplet---unlike the prefactor of the $\log(N)$ term.
Another point to draw attention to is that while the second subcase of \autoref{theo: In/Out Annulus} can be obtained from the first by taking the limit $\eta \to 0^+$ to the first three leading orders, this breaks down at order $\sqrt{N}$, and in fact the constant term diverges in this limit. This is understandable as the expansion is in $N$ not in $\eta$, which is $N$-independent and fixed for a fixed potential $Q$.
\section{Asymptotic expansion of the summands $u_j$}

To set the stage of Laplace's method we rewrite $u_j$ in (\ref{1.18}) according to
\Eq{eq: def uj}{
u_{j}=\int_{0}^{1} re^{-N V_{\tau(j)}(r)}dr,
}
with the parameter $\tau$ depending on $j$,
\begin{equation}\label{eq: def tau}
    \tau=\tau(j):=\frac{j}{N}.
\end{equation}
For the sake of simplicity we will not make the $j$ dependence explicit in the notation. The function $V_{\tau}$ multiplied by $N$ in the exponential term in (\ref{eq: def uj})
is then given by
\begin{equation}\label{3.3}
    V_{\tau}(r):=q(r)-2\tau \log(r).
\end{equation}
To analyse the large $N$ asymptotic of the integral $u_j$ using Laplace's method one needs to determine the location of its minima.

On the open interval $]0,1[$ the minima will be among the solutions of the  critical point equation 
\Eq{eq: SAP}{
    V_{\tau}'(r)=q'(r)-\frac{2\tau}{r}=0. 
}
The potential $Q$ being strictly sub-harmonic  implies that the function $r\mapsto rq'(r)$,
for $0<\tau<1$,
is strictly increasing inside the droplet. Therefore, there exists a unique $r_\tau$ such that
\Eq{}{
 V_{\tau}'(r_\tau)=0,\qquad V_\tau''(r_\tau)>0.
}
The unique critical point $r_\tau$ is then a minimum and satisfies the relation
\Eq{eq: SAP2}{
    r_\tau q'(r_\tau)=2\tau.
}
Note that for $\tau=0$ we have $q'(r_0)=0$. Moreover, the presence of a hard wall is encoded by the existence of $\tau_0\in [0,1[$ such that 
\Eq{eq: def tau0}{
r_{\tau_0}=1,
}
in  keeping with the outer radius $r_1$ satisfying the equation $r_1 q'(r_1) = 2$; recall the text below (\ref{eq: droplet}).

The expression of $r_\tau$ is implicit and so is its derivative with respect to $\tau$. However, one can relate the two by using implicit differentiation on \eqref{eq: SAP2} to arrive at
\Eq{eq: r'}{
    \frac{d r_\tau}{d\tau}=\frac{2}{ r_\tau \Delta Q(r_\tau)}.
}
Note that we have, due to \eqref{eq: SAP2},
\Eq{}{
\Delta Q(r_\tau)=q''(r_\tau) + \frac{2\tau}{r_\tau^2},
}
and in particular,
\Eq{}{
\Delta Q(r_0)=q''(r_0), \quad \Delta Q(1)=q''(1) + q'(1),\quad q'(1)=2\tau_0.
}
Moreover, using \eqref{eq: Laplacian} one can express the successive derivatives of $V_\tau$ in terms of $\Delta Q$
\Eq{eq: derivative relations}{
&V_\tau''(r) = \Delta Q - \frac{V_\tau'(r)}{r}, \quad 
V_\tau^{(3)}(r) = \partial_r \Delta Q - \frac{1}{r} \Delta Q + \frac{2 V_\tau'(r)}{r^2}, \\
& V_\tau^{(4)}(r) = \partial_r^2 \Delta Q + \frac{3}{r^2} \Delta Q - \frac{1}{r} \partial_r \Delta Q - \frac{6 V_\tau'(r)}{r^3}.
}

For the asymptotic expansions of $u_j$ \eqref{eq: def uj}, as $N\to\infty$, there are two different regimes depending on whether the critical point $r_\tau$ \eqref{eq: SAP} is greater than $r_{\tau_0}=1$ or smaller. Recalling the definition of $\tau$ \eqref{eq: def tau} and $\tau_0$ \eqref{eq: def tau0}, for $\tau\leq\tau_0$, we have $r_0\leq r_\tau\leq r_{\tau_0}=1\leq r_1$. For $\tau>\tau_0$ we have $r_0\leq r_{\tau_0}=1<r_\tau\leq r_1$ and thus the critical point lies outside of the integration domain $[0,1]$.

\begin{figure}[H]
        \centering
        \includegraphics[width=0.7\textwidth]{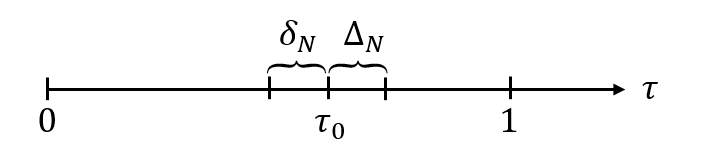}
        \captionsetup{justification=justified}
        \caption{\label{fig: scaling} Illustration of the scalings $\delta_N$ and $\Delta_N$ when $0<\tau_0<1$.}   
\end{figure}

To carry out the analysis we introduce two different scalings for the regimes. The first one, $\Delta_N$, when $\tau_0$ is approached by above, should satisfy
\begin{equation}\label{eq: Delta_N}
    \Delta_N \to 0, \quad \sqrt{N} \Delta_N \to \infty \quad \text{as } N \to \infty.
\end{equation}
The second scaling, $\delta_N$, when $\tau_0$ is approached from below, is specified by
\Eq{eq: delta_N}{
\delta_N := \frac{(\log(N))^2}{\sqrt{N}};
}
see Figure \ref{fig: scaling}.
\vspace*{-0.9cm}

For $\tau\leq \tau_0$, we have $r_\tau\in[r_0,1]$ and one needs to distinguish two additional cases: $r_0 > 0$ (annulus droplet) and $r_0 = 0$ (disk droplet). The analysis is different for $r_0=0$ as the critical point can lie at the lower boundary of the integration domain $[0,1]$.  The details of the latter case are already available in the literature \cite[\S 3.1]{Byun2023}.
\subsection{Case: $\tau\geq\tau_0$}

Recall that for $\tau > \tau_0$ the minimum of
$V_\tau(r)$ in (\ref{eq: def uj}) lies outside the endpoint $r=1$ of the integration interval $[0,1]$. When this separation is beyond a certain threshold as quantified by (\ref{eq: Delta_N}), expansion of the exponent about the endpoint and integration by parts gives an asymptotic expansion in powers of $1/N$. However for distance closer than this threshold, the algebraic portion of the expansion now is in powers of $1/\sqrt{N}$.
The following expansion allows for both these circumstances.

\prop{lem: uj asym0}{
If $\tau\geq\tau_0$, then the integral $u_j$ \eqref{eq: def uj} admits the $N\to\infty$ expansion 
\Eq{eq: uj interpol}{
u_j{=}&\,\frac{ e^{-N V_\tau(1)} }{12 V_\tau''(1)^5 N^2} \Bigg(
-2 V_\tau'(1)^3 V_\tau''(1) V_\tau^{(3)}(1) N^2
+ 2 V_\tau'(1)^2 V_\tau''(1)^2 V_\tau^{(3)}(1) N^2 -10 V_\tau'(1) V_\tau''(1)^2 V_\tau^{(3)}(1) N \\
&
-12 V_\tau''(1)^4 N
    +4 V_\tau''(1)^3 V_\tau^{(3)}(1) N
     + \sqrt{2\pi} \sqrt{V_\tau''(1) N} \, e^{\frac{V_\tau'(1)^2 N}{2 V_\tau''(1)}} 
    \erfc\left(\frac{-V_\tau'(1) \sqrt{N/V_\tau''(1)}}{\sqrt{2}}\right)\\
    &\times\bigg[
        V_\tau'(1)^3 V_\tau^{(3)}(1) N^2 (V_\tau''(1) - V_\tau'(1))+ 6 V_\tau''(1)^3 N (V_\tau''(1) - V_\tau'(1))\\
    &+ 3 V_\tau'(1) V_\tau''(1) V_\tau^{(3)}(1) N (V_\tau''(1) - 2 V_\tau'(1))- 3 V_\tau''(1)^2 V_\tau^{(3)}(1)\bigg]\Bigg)\left(1 + {\rm o}(N c_N^{-4})\right),
}
where $c_N := \sqrt{N} + |V_\tau'(1)|N$.
}

\begin{proof}[Proof of \autoref{lem: uj asym0}]
The first step is to split the integration domain into two subregions 
$
u_j=I_{1}^{\rm out}+I_{2}^{\rm out}
$
where
\begin{equation}
    I_{1}^{\rm out} := \int_{0}^{1 - b_N} r\, e^{-N V_{\tau}(r)}\, dr, \quad 
    I_{2}^{\rm out} := \int_{1 - b_N}^{1} r\, e^{-N V_{\tau}(r)}\, dr,
\end{equation}
with  the cut-off scaling
\begin{equation}
    b_N := \frac{ (\log(N))^2}{c_N}
\end{equation}
and $c_N$ as in \autoref{lem: uj asym0}. The particular choice of the scaling $c_N$ will become clear shortly. Let us first treat the integral $I_{1}^{\rm out}$.

One first notices that $r\mapsto V_\tau(r)$ is strictly decreasing on the interval $[0, 1]$. Thus for $\tau>\tau_0$, the minimum lies at $r=1$. Hence the bound
\begin{equation}
    N(V_\tau(r) - V_\tau(1)) \geq N(V_\tau(1 - b_N) - V_\tau(1)) \geq C_1 (\log(N))^2
\end{equation}
holds for some positive constant \( C_1 > 0 \) and $r\in [0,1-b_N]$. This yields the estimate
\begin{equation}
    I_{1}^{\rm out} \leq e^{-N V_\tau(1)} e^{-C_1 (\log(N))^2},
\end{equation}
telling us that the integral $I_{1}^{\rm out}$ is exponentially suppressed as $N\to\infty$. 

Let us then turn to $I_{2}^{\rm out}$ and proceed with Laplace's method. To this end, we expand $V_\tau(r)$ in a Taylor series around $r=1$ up to third order.

Substituting and proceeding with the change of variable $r - 1 = y/c_N$, the integral becomes
\begin{equation}
\label{eq:integral_change}
    I_{2}^{\rm out}\underset{N\to\infty}{\sim}\frac{e^{-N V_\tau(1)}}{c_N} \int_{-(\log(N))^2}^{0} \left(1+\frac{y}{c_N}\right) e^{-N\left( V_\tau'(1) \frac{y}{c_N} + \frac{V_\tau''(1)}{2} \frac{y^2}{c_N^2} + \frac{V_\tau'''(1)}{6} \frac{y^3}{c_N^3} + {\rm o}(1/N) \right)} dy.
\end{equation}
Note here that the particular choice of $c_N$ is due to $V_\tau'(1)$ vanishing at $\tau=\tau_0$. By choosing $c_N$ in such a way guarantees that either the second or first Taylor coefficient will be exactly of order $1$ for all $\tau\geq\tau_0$, while making the higher order coefficients sub-leading as $N\to\infty$.

As $N y^3/c_N^3\to0$ in the limit $N\to\infty$, we may expand the cubic term to obtain
\Eq{eq:integral_expansion}{
   I_{2}^{\rm out}\underset{N\to\infty}{\sim}\frac{ e^{-N V_\tau(1)} }{c_N}\int_{-(\log(N))^2}^{0}&\, \left(1+ \frac{y}{c_N}\right)
    \left( 1 - \frac{N V_\tau'''(1)}{6} \frac{y^3}{c_N^3} \right)\\
    &\times e^{-N\left( V_\tau'(1) \frac{y}{c_N} + \frac{V_\tau''(1)}{2} \frac{y^2}{c_N^2} \right)}
    \left(1 + {\rm o}(N c_N^{-4})\right) dy.
}
Moreover, since either \( N c_N^{-1} V_\tau'(1) \) or \( N c_N^{-2} V_\tau''(1) \) tends to a positive limit, we can safely extend the domain of integration from \( [-(\log(N))^2, 0] \) to \( (-\infty, 0] \) with exponentially small error. Thus, the integral can be evaluated explicitly and this yields \eqref{eq: uj interpol}.
\end{proof}

With $\tau \ge \tau_0$, we now separate the subcases $\tau - \tau_0 \ge \Delta_N$ and
$\tau - \tau_0 < \Delta_N$ in \autoref{lem: uj asym0}. The general result then simplifies.

\coro{coro: uj asym1}{
If \( \tau - \tau_0 \ge \Delta_N \), then the integral $u_j$ \eqref{eq: def uj} admits the $N\to\infty$ expansion 
\begin{equation}\label{eq: uj asym1}
    u_j = \frac{e^{-N V_\tau(1)}}{N} \left( -\frac{1}{V_\tau'(1)} + \frac{1}{N} \frac{V_\tau''(1) - V_\tau'(1)}{V_\tau'(1)^3} + {\rm o}(1/N) \right).
\end{equation}
}

\begin{proof}[Proof of \autoref{coro: uj asym1}]
 If $ \tau - \tau_0 \ge \Delta_N $, $c_N \gtrsim N \Delta_N$. The order of the error term is then
\Eq{}{
N c_N^{-4} \underset{N\to\infty}{=} {\rm o}(1/N)
}
and the expansion \eqref{eq: uj interpol} reduces to \eqref{eq: uj asym1}.
\end{proof}

\ex{E2}{\normalfont
Set $q(r) = r^2 - 2a \log(r)$, $0<a<1$ in (\ref{3.3}). Then we have
\Eq{3.23}{u_j = \frac{1 }{ 2} N^{-(N(\tau+a) + 1)}
\gamma(N (\tau +a) + 1;N);
}
cf.~(\ref{q3}). For this choice of $q(r)$ the minimum occurs at $r=r_\tau = \sqrt{\tau+a}$,
and thus $\tau_0 = 1-a$. The region $\tau - \tau_0 \ge \Delta_N$ then corresponds to
$\tau + a \ge 1 + \Delta_N$.
Suppose in particular that $\tau + a \ge 1 + \epsilon$ for $\epsilon > 0$ fixed. We then have available a classical asymptotic formula for the incomplete gamma function due to Mahler \cite{Mahler1930} 
(see \cite[Eq.~(2.1)]{Nemes2018} for a contemporary reference) stating that for $z,u \to \infty$ with $u/z > 1$,
\Eq{3.24}{\gamma(u;z) = \frac{z^u e^{-z} }{ u - z}
\bigg ( 1 - \frac{z }{ (u - z)^2} +{\rm O} \Big ( \frac{1 }{ z^2} \Big ) \bigg ). }
Using this in (\ref{3.23}) then shows that under the stated bound on $\tau+a$, for large $N$,
\Eq{3.25}{u_j = \frac{e^{-N} }{ 2 c N} \bigg (
1 - \frac{1 }{ N} \Big ( \frac{1 }{ c^2} + \frac{1 }{ c} \Big ) +
{\rm O} \Big ( \frac{1 }{ N^2} \Big ) \bigg ),
}
where $c:=  \tau + a -1$. On the other hand, using that
\Eq{}{
V_\tau(1), \quad V_\tau'(1) = 2c, \quad V_\tau''(1) = - 2 c + 4,
}
the result of \autoref{coro: uj asym1} also gives an explicit asymptotic formula for the first two terms of $u_j$ in this regime. Precise agreement with (\ref{3.25}) is displayed.
}
\coro{coro: uj asym2}{
If \( \tau - \tau_0 < \Delta_N \), then the integral $u_j$ \eqref{eq: def uj} admits the $N\to\infty$  expansion
\begin{equation}\label{eq: uj asym2}
     u_j = \frac{e^{-N V_\tau(1)}}{N} \sqrt{ \frac{2 N}{\Delta Q(1)} } \left( f_1(x) + \frac{1}{\sqrt{N}}\sqrt{ \frac{2 }{\Delta Q(1)} }  h_1(x) + {\rm O}(N^{-1}) \right),
\end{equation}
where 
\Eq{}{
x := \sqrt{\frac{2 N }{\Delta Q(1)}}(\tau - \tau_0) \geq 0,
}
and the functions $f_1$ and $h_1$ are given by
\Eq{}{
f_1(x):=\frac{\sqrt{\pi}}{2}e^{x^2}\erfc(x),
}
\Eq{3.30}{
h_1(x):=\frac{1}{3}\left[(x^2-2)-x(2x^2-3)f_1(x) \right]+\frac{1}{6}\frac{\partial_r\Delta Q(1)}{\Delta Q(1)}\left[(x^2+1)-x(2x^2+3)f_1(x) \right].
}
}

\begin{proof}[Proof of \autoref{coro: uj asym2}]
If $ \tau - \tau_0 < \Delta_N $. The order of the error term is then
\Eq{}{
N c_N^{-4} \underset{N\to\infty}{=} {\rm O}(1/N),
}
and the expansion \eqref{eq: uj interpol} reduces to \eqref{eq: uj asym2}.
\end{proof}

\ex{E3}{\normalfont ${}$
\\
1.~Let us choose $q(r) = r^2 - 2 a \log(r)$, $0<a<1$ as in 
Example \ref{E2}. The significance of this is that we then have available the special function evaluation (\ref{3.23}). The requirement $\tau - \tau_0 < \Delta_N$ for this $q(r)$ then reads $\tau + a < 1 + d_0/N^{1/2}$ for some $d_0 > 0$. Relevant to this parameter regime is the incomplete gamma function $y \to \infty$ asymptotic formula \cite{Tricomi1950},
\cite[\S 2]{Nemes2018}
\Eq{}{
\frac{\gamma(y+1; y - t y^{1/2}) }{ \Gamma(y+1)}
=\frac{1 }{ 2}\erfc(2^{-1/2} t) - \frac{1 }{ \sqrt{2 \pi y} } \exp \Big ( - \frac{t^2 }{ 2} \Big )\frac{t^2 + 2 }{ 3} + {\rm O} (y^{-1}),
}
uniformly valid for bounded $t$. Minor manipulation of this shows  that for $N \to \infty$
\Eq{3.33}{
 \frac{\gamma(N + \sqrt{2N} x + 1;N) }{
 \Gamma(N + \sqrt{2N} x + 1) }
 = \frac{1 }{ 2} \erfc(x) + \frac{1 }{ 3} \sqrt{\frac{1 }{ 2 \pi N}} e^{-x^2} (x^2 - 2) +
{\rm O}(N^{-1}),
}
which with the choice
$x = (N/2)^{1/2}(\tau - \tau_0)$
is directly relevant to (\ref{3.23}). Combining this with Stirling's formula we deduce the large $N$ expansion
\Eq{}{u_j=  \sqrt\frac{\pi }{ 2  N}
\frac{e^{-N(\tau + a)} }{ 2} \bigg (
\frac{1 }{ 2} \erfc(x) + 
\frac{1 }{ 2} \sqrt\frac{1 }{ 2 N} \erfc(x)
\Big ( x - \frac{2 x^3 }{ 3} \Big )
+ \frac{1 }{ 3} \sqrt\frac{1 }{ 2 \pi N} e^{-x^2} (x^2 - 2) +{\rm O}(N^{-1}) \bigg ).
}
This is in precise agreement with the appropriate specialisation of \autoref{coro: uj asym2}. \\
2.~A drawback with the above example is that $\partial_r \Delta Q(r) = 0$ so the second main term in (\ref{3.30}) vanishes. A choice of radial potential where this is not the case is
$q(r) = r^{2(1+\mu)} - 2a \log(r)$, $a > 0$ (the case $\mu = 0$ is that of 1.~above). Now $\tau_0 =
1 + \mu - a < 1$ and it is required that $0 < \tau_0 < 1$.
This choice remains in the class that permits an evaluation in terms of the incomplete gamma function for $u_j$,
\Eq{}{
 u_j = \frac{1 }{ 2 (1 + \mu)}
 N^{-(1+ N(a+\tau))/(1+\mu)}
 \gamma((-\mu+ N(a+\tau))/(1+\mu)+1;N).}
 Writing $N(a+\tau)/(1+\mu) = N +
 \sqrt{2N} x$ with $x = (1 + \mu)^{-1} (N/2)^{1/2} (\tau - \tau_0)$,
 we can again make use of (\ref{3.33}) to reclaim the large $N$ expansion of
 \autoref{coro: uj asym2}. Moreover, extensive use was made of this example by way of numerical evaluations associated with the first subcase of
\autoref{theo: In/Out Annulus},  to confirm the accuracy of various terms as stated.
}

\subsection{Case: $\tau\leq\tau_0$ ($r_0>0$ annulus support)}
As for the case $\tau \ge \tau_0$, the large $N$ expansion depends on $|\tau - \tau_0|$, although now with the quantity $\delta_N$ defined in (\ref{eq: delta_N}) being the length scale determining a transition in behaviours. We consider first the case $|\tau - \tau_0|< \delta_N$, with our result expressed in terms of 
the Taylor expansion of
$V_\tau(r)$
about $r= r_\tau$ up to fourth order.
Once this is established, we present the same expansion, now using $\Delta Q(r)$ and its derivatives at $r=1$, as a corollary.
\prop{lem: uj asym3}{
In the case $r_0>0$ and for $\tau\leq\tau_0$ and $\abs{\tau - \tau_0} \leq \delta_N$, the integral $u_j$ \eqref{eq: def uj} admits the $N\to\infty$ expansion
\Eq{eq: uj interpol2}{
&u_j {=}\frac{\exp(-N V_\tau(r_\tau)-\frac{1}{2} N (r_\tau-1)^2 V_\tau''(r_\tau))}{144 N^{3/2} V_\tau''(r_\tau)^3} \Bigg[9 \sqrt{2 \pi } \sqrt{V_\tau''(r_\tau)} \exp(\frac{1}{2} N (r_\tau-1)^2 V_\tau''(r_\tau))\\
&\times\left(8 N r_\tau V_\tau''(r_\tau)^2-r_\tau V_\tau^{(4)}(r_\tau)-4 V_\tau^{(3)}(r_\tau)\right)\erfc\left(\frac{(r_\tau-1) \sqrt{N V_\tau''(r_\tau)}}{\sqrt{2}}\right)\\
&+2 \sqrt{N} \bigg[\bigg(N (r_\tau-1)^2 \bigg[V_\tau^{(3)}(r_\tau) \left(N (r_\tau-1)^3 V_\tau^{(3)}(r_\tau)+12\right)-3 (r_\tau-1) V_\tau^{(4)}(r_\tau)\bigg]-72\bigg) V_\tau''(r_\tau)^2 \\
&+\bigg(-N (r_\tau-6) (r_\tau-1)^3 V_\tau^{(3)}(r_\tau)^2
+3 (r_\tau-4) (r_\tau-1) V_\tau^{(4)}(r_\tau)-12 (r_\tau-3) V_\tau^{(3)}(r_\tau)\bigg) V_\tau''(r_\tau)\\
&-3 (r_\tau-1) (3 r_\tau-8) V_\tau^{(3)}(r_\tau)^2\bigg]+{\rm O}\left(N^{-1/2} \right)\Bigg].
}
}

\begin{proof}[Proof of \autoref{lem: uj asym3}]
Here $\tau \leq \tau_0$ implies that the critical point, which is also the unique minimum, lies inside the domain of integration, i.e. $r_\tau \in[r_0,1]\subset[0,1]$, with $r_0>0$. Assuming $\abs{\tau - \tau_0} \leq \delta_N$, where $\delta_N$ is given by \eqref{eq: delta_N}, we decompose the integral $u_j$ as
\Eq{}{
u_j=I_{1}^{\rm in}+I_{2}^{\rm in},
}
where
\begin{equation}\label{eq: I1 I2}
    I_{1}^{\rm in} := \int_{0}^{r_\tau - \beta_N} r\, e^{-N V_{\tau}(r)}\, dr, \quad 
    I_{2}^{\rm in} := \int_{r_\tau - \beta_N}^{1} r\, e^{-N V_{\tau}(r)}\, dr,
\end{equation}
with the cut-off scaling defined by
\Eq{eq: beta_N}{
\beta_N := \frac{\log(N)}{\sqrt{N}}.
}

Let us start with $ I_{1}^{\rm in}$. As the critical point $r_\tau$ is outside the interval $[0,r_\tau-\beta_N]$, the integrand has a unique minimum attained at $r=r_\tau-\beta_N$. As a consequence,  $r\mapsto V_\tau(r)$ is strictly decreasing on $[0,r_\tau-\beta_N]$. Therefore, one has
\begin{equation}\label{eq: bound I1}
    N\left(V_\tau(r) - V_\tau(r_\tau)\right)
    \geq N\left(V_\tau(r_\tau - \beta_N) - V_\tau(r_\tau)\right)
    \geq C (\log(N))^2,
\end{equation}
for some constant $C > 0$. This implies the bound
\begin{equation}
     I_{1}^{\rm in}
    \leq e^{-N V_\tau(r_\tau)} \, e^{-C (\log(N))^2},
    \end{equation}
and so the integral $I_{1}^{\rm in}$ is exponentially suppressed as $N\to\infty$.

Let us turn to consider $I_{2}^{\rm in}$ and proceed with Laplace's method. To this end, we expand $V_\tau(r)$ in a Taylor series around $r_\tau$
to fourth order.

Proceeding with the change of variable $r = r_\tau + \frac{t}{\sqrt{N}}$ and then Taylor expanding the exponential yields
\begin{equation}
\begin{aligned}
    e^{-N V_\tau(r)}
    \underset{N\to\infty}{=}&\, e^{-N V_\tau(r_\tau)} \, e^{-V_\tau''(r_\tau) t^2}\\
    &\times\Bigg(
        1 - \frac{1}{6} \sqrt{\frac{1}{N}} V_\tau^{(3)}(r_\tau) t^3
        + \frac{V_\tau^{(3)}(r_\tau)^2 t^6 - 3 V_\tau^{(4)}(r_\tau) t^4}{72 N}+{\rm o}(N^{-1})
    \Bigg).
\end{aligned}
\end{equation}
Substituting in the integrand for $I_{2}^{\rm in}$ implies
\begin{equation}
\begin{aligned}
    I_{2}^{\rm in}
    \underset{N\to\infty}{=}&\,  \frac{e^{-N V_\tau(r_\tau)}}{\sqrt{N}} \int_{-\log(N)}^{(1 - r_\tau) \sqrt{N}}
    \left(r_\tau + \frac{t}{\sqrt{N}}\right) e^{-V_\tau''(r_\tau) t^2}\\
    & \times \left(
        1 - \frac{1}{6} \sqrt{\frac{1}{N}} V_\tau^{(3)}(r_\tau) t^3
        + \frac{V_\tau^{(3)}(r_\tau)^2 t^6 - 3 V_\tau^{(4)}(r_\tau) t^4}{72 N}+{\rm O}(N^{-3/2})
    \right)
     \, dt.
\end{aligned}
\end{equation}
Since $V_\tau''(r_\tau)>0$ we can extend the integration domain from $[-\log(N),0]$, to $(-\infty,0]$ with an exponentially small error. The integral can be computed explicitly and gives the claim \eqref{eq: uj interpol2}.
\end{proof}

\coro{coro: uj asym4}{
In the case $r_0>0$ and for $\tau\leq\tau_0$ and $\abs{\tau - \tau_0} \leq \delta_N$, the integral $u_j$ \eqref{eq: def uj} admits the $N\to\infty$ expansion
\begin{equation}\label{eq: uj asym4}
     u_j =  \frac{\sqrt{2\pi}r_\tau}{\sqrt{N V_\tau''(r_\tau)}}e^{-N V_\tau(r_\tau)}\left(f_2(x)+ \frac{1}{\sqrt{N}}\sqrt{ \frac{2 }{\Delta Q(1)} }h_2(x)+{\rm O}(N^{-1})\right),
\end{equation}
where 
\Eq{eq: def x}{
x := \sqrt{\frac{2 N }{\Delta Q(1)}}(\tau - \tau_0) \leq 0,
}
and the functions $f_2$ and $h_2$ are given by
\Eq{}{
f_2(x):=\frac{1}{2}\erfc(x), \quad
h_2(x):=\frac{e^{-x^2} \left(x^2-2\right)}{3 \sqrt{\pi }}+\frac{e^{-x^2} \left(x^2+1\right) }{6 \sqrt{\pi } }\frac{\partial_r \Delta Q(1)}{\Delta Q(1)}.
}
}

\begin{proof}[Proof of \autoref{coro: uj asym4}]
As $\abs{\tau - \tau_0} \leq \delta_N$, with $\delta_N$ given in \eqref{eq: delta_N}, one has the following expansions, as $\abs{\tau - \tau_0}\to0$ and thus as $N\to\infty$,
\Eq{eq: expansions V r}{
V_\tau''(r_\tau)&=\Delta Q(1)+2\frac{\partial_r \Delta Q(1)}{\Delta Q(1)} (\tau-\tau_0)+{\rm O}(\tau-\tau_0)^2,\\
r_\tau &=1+\frac{2}{\Delta Q(1)}(\tau-\tau_0)-2\frac{ \Delta Q(1)+\partial_r \Delta Q(1)}{\Delta Q(1)^3} (\tau-\tau_0)^2+{\rm O}(\tau-\tau_0)^3.
}
Then, using the new variable $x$ defined in \eqref{eq: def x}, the expansion \eqref{eq: uj interpol2} reduces to \eqref{eq: uj asym4}.
\end{proof}

\rem{} 
{\normalfont For purposes of matching with the $\tau \le \tau_0$, $|\tau - \tau_0| \ge \delta_N$ case to be considered next, the leading factor in 
(\ref{eq: uj asym4}) has been written in terms of $r_\tau, V_\tau(r_\tau)$, $V_\tau''(r_\tau)$,
rather than the expansion of these terms about $r_{\tau_0} = 1$. If one was to carry out the expansions to order $1/\sqrt{N}$ as implied by 
(\ref{eq: expansions V r}), the asymptotic expansion of  \autoref{coro: uj asym4} becomes identical to that of  \autoref{coro: uj asym2}.
}

\medskip

As said, we proceed to consider the circumstance that $\abs{\tau - \tau_0} \geq \delta_N$. Actually, the required working has already been carried out in \cite[Lemma 3.2]{Byun2023}, where in particular the quantity $\mathcal B(r)$ in (\ref{eq: def B}) below was first identified.

\prop{lem: uj asym5}{
In the case $r_0>0$ and for $\tau\leq\tau_0$, if $\abs{\tau - \tau_0} \geq \delta_N$ , then the integral $u_j$ \eqref{eq: def uj} admits the $N\to\infty$ expansion
\begin{equation}\label{eq: uj asym5}
     u_j\underset{N\to\infty}{=}\frac{\sqrt{2\pi}r_\tau}{\sqrt{N V_\tau''(r_\tau)}}e^{-NV_\tau(r_\tau)}\left(1+\frac{1}{N} \mathcal{B}(r_\tau)+{\rm O}(N^{-2})\right),
\end{equation}
where 
\Eq{eq: def B}{
\mathcal{B}(r):=- \frac{1}{8} \frac{\partial_r^2 \Delta Q(r)}{\left(\Delta Q(r)\right)}
- \frac{19}{24 r} \frac{\partial_r \Delta Q(r)}{\left(\Delta Q(r)\right)^2}
+ \frac{5}{24} \frac{\left(\partial_r \Delta Q(r)\right)^2}{\left(\Delta Q(r)\right)^3}
+ \frac{1}{3r^2} \frac{1}{\Delta Q(r)}.
}
}

\begin{proof}[Proof of \autoref{lem: uj asym5}] (Cf.~\cite[proof of Lemma 3.2]{Byun2023})
Using the fact $\abs{\tau - \tau_0} \geq \delta_N$ , we decompose the integral defining $u_j$ according to
\Eq{3.50}{
u_j=I_{1}^{\rm in}+I_{2-}^{\rm in}+I_{2+}^{\rm in},
}
where $I_{1}^{\rm in}$ is given by \eqref{eq: I1 I2} and
\Eq{}{
    I_{2-}^{\rm in} := \int_{r_\tau - \beta_N}^{r_\tau + \beta_N} r\, e^{-N V_{\tau}(r)}\, dr, \quad 
    I_{2+}^{\rm in} := \int_{r_\tau + \beta_N}^{1} r\, e^{-N V_{\tau}(r)}\, dr,\\
}
with $\beta_N$ is given by \eqref{eq: beta_N}.

We note that $r_\tau + \beta_N < 1$. Indeed, $\tau\mapsto r_\tau$ is strictly increasing due to \eqref{eq: r'} and the sub-harmonicity of $Q$. Therefore, for $\tau\leq\tau_0$  and $\abs{\tau - \tau_0} \geq \delta_N$, we have,
\Eq{}{
1-r_\tau\geq 1-r_{\tau_0-\delta_N}.
}
Thus, from \eqref{eq: expansions V r} for $N$ large enough, 
\Eq{}{
1-r_\tau\geq 1-r_{\tau_0-\delta_N}\geq \frac{2}{\Delta Q(1)} \delta_N>\beta_N>0.
}

We consider now the $N\to\infty$ form of the integrals (\ref{3.50}). 
The integral $I_{1}^{\rm in}$ remains exponentially suppressed as $N\to\infty$, since \eqref{eq: bound I1} still holds here. One can then move to $I_{2+}^{\rm in}$. A Taylor expansion of $r\mapsto V_\tau(r)$ around $r=r_\tau$ yields the inequalities
\begin{equation}
    N\left(V_\tau(r) - V_\tau(r_\tau)\right) 
    \geq N V_\tau''(r_\tau) b_N^2 
    \geq C_1 \log^2 N,
\end{equation}
where $C_1$ is some positive constant. Consequently, we have the bound
\begin{equation}
    I_{2+}^{\rm in}
    \leq e^{-N V_\tau(r_\tau)} e^{-C_1 \log^2 N}.
\end{equation}
Thus $ I_{2+}^{\rm in}$ is also exponentially suppressed as $N\to\infty$. The only non-negligible contribution then comes from the integral $ I_{2-}^{\rm in}$. The large $N$ expansion for this is obtained via the classic Laplace's method \cite{Wong2001} and yields \eqref{eq: uj asym5}. 

\end{proof}

\ex{}{ \normalfont 
Consider the choice of $q(r)$ in Example \ref{E2}, for which $u_j$ is given in terms of the incomplete gamma function according to (\ref{3.23}). Relevant to the regime of 
 \autoref{lem: uj asym5} is the asymptotic formula \cite[equivalent to Eq.~8.11.7]{DLMF},
giving that for $z,u \to \infty$ with $u/z < 1$,
\Eq{}{\gamma(u;z) =\Gamma(u)+ \frac{z^u e^{-z} }{ u - z}
\bigg ( 1 - \frac{z }{ (u - z)^2} +{\rm O} \Big ( \frac{1 }{ z^2} \Big ) \bigg ); }
cf.~(\ref{3.24}). In fact the second term on the right hand side is exponentially suppressed relative to $\Gamma(u)$, so in this regime the asymptotics follow from the logarithm of Stirling's formula, extended to include the order $1/U$ term. Noting that $r_\tau = (a+\tau)^{1/2}$, $V_\tau''(r_\tau) = \Delta Q(r_\tau) = 4$,
$\partial_r \Delta Q(r) = 0$, agreement is found with \autoref{lem: uj asym5}.
}

\subsection{Case: $\tau\leq\tau_0$ ($r_0=0$ disk support)}\label{S3.3}

As mentioned previously, in the case of a disk droplet ($r_0 = 0$) the critical point $r_\tau$ can be at the lower boundary of the domain. This requires an additional local analysis which can be found in \cite[\S 3.1]{Byun2023}. Here again, one needs to control how close $r_\tau$ is from $0$ depending on a $N$-dependent scaling of the paramater $\tau=j/N$. 

In \cite[\S 3.1]{Byun2023} this scale is chosen to be of order ${\rm O}(N^{\epsilon - 1})$ with $0 < \epsilon < 1/5$. This yields the following asymptotic expansion for $u_j$.

\prop{lem: uj asym6}{(\cite[Lemma 3.1, with $h_j = 2 u_j$]{Byun2023})
In the case $r_0=0$, for $\tau \le {\rm O}(N^{\epsilon - 1})$, $(0 < \epsilon < 1/5$) the integral $u_j$ \eqref{eq: def uj} admits the $N\to\infty$ expansion
\begin{equation}\label{eq: uj asym6}
     2 u_j\underset{N\to\infty}{=} -2 N q(0) - (j+1) \log (N q''(0)) + \log(j!) +{\rm O}(N  \tau^{3/2}(\log(N))^3).
\end{equation}
}

Although we refer to the original reference for the details of the proof, we remark that it is a classical Laplace's method applied when the critical point is exactly at the boundary of the domain. This yields only half of the Gaussian integral approximation compared to the case when the critical point lies inside the domain. 

\section{Proof In/Out Annulus case: $0 < r_0 < 1 \leq r_1$}

As seen in the previous section, the asymptotic behaviour of the terms $u_j$ \eqref{eq: def uj} depends on the range of the indices $j$. Specifically, the asymptotic expansion is obtained via Laplace's method and one has that for indices larger than $N\tau_0$, the critical point escapes the domain of integration. The global maximum of the integrand is then located at the right boundary of the integration domain. This justifies the splitting of the sum
\Eq{4.1c}{
\log(\frac{{Z}_{N}^{\rm h}[q(r)]}{(2 \pi)^N })=\sum_{j = 0}^{\lfloor N\tau_0 \rfloor - 1} \log u_j+\sum_{j = \lfloor N\tau_0 \rfloor}^{N-1} \log u_j,
}
where $\lfloor .\rfloor$ is the floor function. In addition to the asymptotic formulas for $u_j$ in the various regimes of the previous section, crucial at this stage is the Euler-Maclaurin summation formula. This is revised in \autoref{appendix}. It is used to obtain the large $N$ expansion of the sums, up to terms which go to zero in this limit.

\prop{lem: sum out}{
For $\tau_0<1$ we have the  $N\to\infty$ expansion
\Eq{eq: sum out}{
\sum_{j = \lfloor N\tau_0 \rfloor }^{N-1} \log u_j =&\,-(\eta N - \{N \tau_0\} )(N q(1) +  \log(N)) -N \eta (\log (2 \eta)-1)-\sqrt{N}\gamma^{\rm out}  \sqrt{\Delta Q(1) }-\frac{\log(N)}{4}\\
&-\alpha^{\rm out} +\beta^{\rm out}\frac{  \partial_r\Delta Q(1)}{\Delta Q(1)}+\frac{\Delta Q(1)}{4 \eta }-\frac{1}{2} \log(\eta)+\frac{1}{4} \log(\frac{\pi  \Delta Q(1)}{2})+{\rm o}(1),\\
}
where the universal positive constants $\alpha^{\rm out}, \beta^{\rm out}, \gamma^{\rm out}$ are given in \eqref{eq: pos constant out}, and where we denote the fractional part $\{x\} :=
x - \lfloor x \rfloor$, for $x\geq 0$.
}

\begin{proof}[Proof of \autoref{lem: sum out}]
Let us start by splitting the sum 
\Eq{}{
S^{\rm out}:=\sum_{j = \lfloor N\tau_0 \rfloor}^{N-1} \log u_j =&\,\sum_{j = \lfloor N\tau_0 \rfloor }^{\lfloor N(\tau_0 + \Delta_N) \rfloor} \log u_j +\sum_{j = \lfloor N(\tau_0 + \Delta_N) \rfloor +1}^{N-1}\log  u_j.
}
To exploit the asymptotics of $u_j$ given in \autoref{coro: uj asym1} and \autoref{coro: uj asym2}, we subdivide the sum further by writing
\Eq{}{
S^{\rm out}=S^{\rm out}_0+S^{\rm out}_{1-}+S^{\rm out}_{2-}+S^{\rm out}_{1+}+S^{\rm out}_{2+},
}
where
\Eq{}{
&S^{\rm out}_0:=\sum_{j = \lfloor N\tau_0 \rfloor }^{N} \log(\frac{\exp(-N V_\tau(1))}{N}),\quad S^{\rm out}_{1-}:=\sum_{j = \lfloor N\tau_0 \rfloor }^{\lfloor N(\tau_0 + \Delta_N) \rfloor }\log(\sqrt{ \frac{2 N}{\Delta Q(1)} } f_1(x)),\\
&S^{\rm out}_{2-}:=\sum_{j = \lfloor N\tau_0 \rfloor }^{\lfloor N(\tau_0 + \Delta_N) \rfloor }\log(1 + \sqrt{ \frac{2 }{N\Delta Q(1)} }  \frac{h_1(x)}{f_1(x)}),\quad S^{\rm out}_{1+}:=\sum_{j = \lfloor N(\tau_0 + \Delta_N) \rfloor +1}^{N-1} \log(\frac{-1}{V_\tau'(1)}),\\
& S^{\rm out}_{2+}:=\sum_{j = \lfloor N(\tau_0 + \Delta_N) \rfloor +1}^{N-1} \log(1+\frac{1}{N}  \frac{V_\tau'(1) - V_\tau''(1)}{(V_\tau'(1))^2}),
}
with $x$ defined as in \eqref{eq: def x}.
\subsubsection*{Computation of $S^{\rm out}_0$}
Let us start by computing $S^{\rm out}_0$. Using the fact $\tau_0=1-\eta$, and $V_\tau(1)=q(1)$, it is immediate that
\Eq{eq: S0}{
S^{\rm out}_0=\sum_{j = \lfloor N\tau_0 \rfloor}^{N-1} \left(-N V_\tau(1) - \log(N)\right)=&\, (-N \eta - \{N \tau_0\})
(N  q(1) +  \log(N)),
}
where $\{.\}$ denotes the fractional part as defined in  \autoref{lem: sum out}.
\subsubsection*{Computation of $S^{\rm out}_{1+}$}
Let us set $y(j):=\tau-\tau_0=j/N-\tau_0$. We thus have $-V_\tau'(1)=2y(j)$ and hence
\Eq{}{
S^{\rm out}_{1+}=\sum_{j = \lfloor N(\tau_0 + \Delta_N) \rfloor +1}^{N-1} -\log(-V_\tau'(1))=-\sum_{j = \lfloor N(\tau_0 + \Delta_N) \rfloor +1}^{N-1}\log(2y(j)).
}
Applying the Euler–Maclaurin formula \autoref{prop: EM formula}, we get
\Eq{4.8}{
  S^{\rm out}_{1+}\underset{N\to\infty}{=}\,  -N \int_{\xi_N^*/N}^{\eta} \log(2y) \, dy 
  + \frac{1}{2}\left(\log(2\eta)+\log(2 \xi_N^*/N)  \right) + {\rm o}(1),
  }
  where we introduce
\Eq{eq: def xi_N}{
\xi_N^*:=&\,\lfloor N(\tau_0 + \Delta_N) \rfloor-N\tau_0= N \Delta_N- \{N(\tau_0 + \Delta_N)\}.
}
Evaluating the integral and simplifying yields
\begin{equation}\label{eq: S1+}
     S^{\rm out}_{1+}\underset{N\to\infty}{=}\,-N \eta\left(\log(2\eta)-1 \right) 
    + \xi_N^*\left( \log(\frac{2\xi_N^*}{N})-1\right) \\
    + \frac{1}{2}\left(\log(2\eta)+\log(\frac{2\xi_N^*}{N})  \right) + {\rm o}(1).
\end{equation}
\subsubsection*{Computation of $S^{\rm out}_{2+}$}
Keeping $y=y(j)$ as defined previously, we have $V_\tau''(1) = \Delta Q(1) + 2y$ and expanding the $\log$ yields
\begin{equation}
    S^{\rm out}_{2+}\underset{N\to\infty}{=} \sum_{j = \lfloor N(\tau_0  + \Delta_N)\rfloor+1}^{N-1} \frac{1}{N} \frac{-\Delta Q(1) - 4y(j)}{4y(j)^2}+{\rm O}(N^{-1}).
\end{equation}
We then define the function
\Eq{}{
g(y(j)) = \frac{-\Delta Q(1) - 4y(j)}{4y(j)^2},
}
and use the Euler–Maclaurin formula \autoref{prop: EM formula} once again to get
\begin{equation}
    S^{\rm out}_{2+}\underset{N\to\infty}{=}\frac{1}{N} \int_{\lfloor N(\tau_0  + \Delta_N)\rfloor}^{N} g(y(s)) \, ds - \frac{1}{2N} \big( g(y(N)) + g(y(\lfloor N(\tau_0  + \Delta_N)\rfloor)) \big) + {\rm O}(N^{-1}).
\end{equation}
The boundary terms can be neglected. Indeed, by assumption, $\sqrt{N}\Delta_N\to0$ \eqref{eq: Delta_N}, thence
\Eq{}{
\frac{1}{2N} \big( g(y(N)) + g(y(\lfloor N(\tau_0  + \Delta_N)\rfloor)) \big)\underset{N\to\infty}{=} {\rm O}\left(\frac{1}{N\Delta_N^2} \right) \underset{N\to\infty}{=} {\rm o}(1).
}
Then, proceeding with the change of variables $s \mapsto s/N -\tau_0 =:Y$,  we obtain
\begin{equation}
    \frac{1}{N} \int_{\lfloor N(\tau_0  + \Delta_N)\rfloor}^{N} g(y(s)) \, ds 
    = \int_{\xi_N^*/N}^{\eta} g(Y) \, dY \\
    = -\frac{\Delta Q(1)}{4} \left( \frac{N}{\xi_N^*} - \frac{1}{\eta} \right) 
    + \log(\frac{\xi_N^*}{N}) - \log(\eta).
\end{equation}
Thus 
 \Eq{eq: S2+}{
 S^{\rm out}_{2+}\underset{N\to\infty}{=} -\frac{\Delta Q(1)}{4} \left( \frac{N}{\xi_N^*} - \frac{1}{\eta} \right) 
    + \log(\frac{\xi_N^*}{N}) - \log(\eta) +{\rm o}(1).
 }

\subsubsection*{Computation of $S^{\rm out}_{1-}$}
From the additive property of the $\log$ of a product, one can rewrite
\Eq{}{
S^{\rm out}_{1-}=\sum_{j = \lfloor N\tau_0 \rfloor}^{\lfloor N(\tau_0 + \Delta_N) \rfloor} \log(\sqrt{ \frac{2 N}{\Delta Q(1)} } )+ \sum_{ j = \lfloor N\tau_0 \rfloor }^{\lfloor N(\tau_0 + \Delta_N) \rfloor} \log(f_1(x(j))).
}
The first sum is immediate,
\Eq{eq: sum1}{
\sum_{j = \lfloor N\tau_0 \rfloor}^{\lfloor N(\tau_0 + \Delta_N) \rfloor } \log(\sqrt{ \frac{2 N}{\Delta Q(1)} } )=(\xi_N +1 )\log(\sqrt{ \frac{2 N}{\Delta Q(1)} } ),
}
where $
\xi_N:=\,\xi_N^*+ \{N\tau_0 \}=\lfloor N(\tau_0 + \Delta_N) \rfloor-\lfloor N\tau_0\rfloor$.
For the second sum, we apply the Euler–Maclaurin formula \autoref{prop: EM formula} and proceed with the change of variable
\Eq{}{
j\mapsto x(j)= \sqrt{\frac{2  }{N\Delta Q(1)}}(j -\lfloor N \tau_0\rfloor).
}
This gives
\begin{equation}\label{eq: sum2}
\begin{aligned}
    \sum_{j = \lfloor N\tau_0 \rfloor}^{\lfloor N(\tau_0 + \Delta_N) \rfloor} \log(f_1(x(j))) 
    &\underset{N\to\infty}{=} \sqrt{\frac{N\Delta Q(1)}{2 }} \int_0^{\sqrt{\frac{2 }{N \Delta Q(1)}}\xi_N} \left( \log(2x f_1(x)) - \log(2x) \right) dx \\
    &\quad + \frac{1}{2} \left[ \log(f_1(0)) + \log(f_1\left(\sqrt{\frac{2  }{N\Delta Q(1)}} \xi_N\right)) \right] + {\rm o}(1),
\end{aligned}
\end{equation}
where 
on the right hand side we have added and subtracted $\log(2x)$ in the integrand to get the function $x\mapsto\log(2x f_1(x))$ which is integrable on $\mathbb{R}_+$. Indeed, using the asymptotic expansion of the complementary error function
\Eq{eq: error fct asymp}{
\erfc(x)& \sim {\frac {e^{-x^{2}}}{x{\sqrt {\pi }}}}\left(1+\sum _{k=1}^{\infty }(-1)^{k}{\frac { (2k-1)!!}{\left(2x^{2}\right)^{k}}}\right),\quad x\rightarrow \infty,
}
we have 
\Eq{4.21}{
\log(2x f_1(x))=\log(\sqrt{\pi}xe^{x^2}\erfc(x))\underset{x\to\infty}{=}-\frac{1}{2x^2}+{\rm O}(x^{-4}).
}

Focusing on the integral on the right hand side of (\ref{eq: sum2}), one can split the two terms. For the first term, one can extend the integration domain to $\infty$ and estimate the error with the help of (\ref{4.21}). Thus we have
\Eq{eq: sum3}{
&\sqrt{\frac{N\Delta Q(1)}{2 }} \int_0^{\sqrt{\frac{2 }{N \Delta Q(1)}}\xi_N}\log(2x f_1(x)) \, dx\\
= &\,\sqrt{\frac{N\Delta Q(1)}{2 }} \left(\int_0^\infty-\int_{\sqrt{\frac{2 }{N \Delta Q(1)}}\xi_N}^\infty \right) \log(2x f_1(x)) \, dx\\
\underset{N\to\infty}{=}&\,\sqrt{\frac{N\Delta Q(1)}{2 }} \int_0^\infty \log(2x f_1(x)) \, dx+ \frac{\Delta Q(1)}{4}\frac{N}{\xi_N} + {\rm o}(1).
}
The remaining portion of the integral on the right hand side of (\ref{eq: sum2})
has the exact evaluation
\Eq{eq: sum4}{
 -\sqrt{\frac{N\Delta Q(1)}{2 }} \int_0^{\sqrt{\frac{2 }{N\Delta Q(1)}}\xi_N} \log(2x) \,dx= \xi_N\left[1-\log(2\sqrt{\frac{2 }{N \Delta Q(1)}}\xi_N)\right]. \\
}
Considering now the $N$-dependent boundary term on the right hand side of (\ref{eq: sum2}), one has 
\Eq{eq: sum5}{
\log(f_1\left(\sqrt{\frac{2  }{N\Delta Q(1)}} \xi_N\right))  \underset{N\to\infty}{=} -\log(2\sqrt{\frac{2  }{N\Delta Q(1)}} \xi_N) + {\rm o}(1).
}

Finally, gathering \eqref{eq: sum1} and \eqref{eq: sum2} with the expressions \eqref{eq: sum3}, \eqref{eq: sum4} and \eqref{eq: sum5} yields
\Eq{}{
S^{\rm out}_{1-}\underset{N\to\infty}{=}&\,(\xi_N +1 )\log(\sqrt{ \frac{2 N}{\Delta Q(1)}})-\sqrt{N\Delta Q(1)} \gamma^{\rm out}\\
&+ \frac{\Delta Q(1)}{4}\frac{N}{\xi_N}+ \xi_N\left[1-\log(2\sqrt{\frac{2 }{N \Delta Q(1)}}\xi_N)\right]
+ \frac{1}{2} \left[ \log(\frac{\sqrt{\pi}}{2}) -\log(2\sqrt{\frac{2  }{N\Delta Q(1)}} \xi_N) \right] + {\rm o}(1),
}
where we have used the definition of the constant $\gamma^{\rm out}$ \eqref{eq: pos constant out}. Simplifying further we arrive at
\Eq{eq: S1-}{
S^{\rm out}_{1-}\underset{N\to\infty}{=}&\,-\gamma^{\rm out}  \sqrt{\Delta Q(1) N}+\frac{1}{4} \log(N)-\frac{1}{4} \log \left(\frac{2 \Delta Q(1)}{\pi }\right)\\
&- \xi_N \left(\log(\frac{2\xi_N}{N})-1\right)+\frac{\Delta Q(1)}{4}\frac{N}{\xi_N}-\frac{1}{2} \log(\frac{2\xi_N}{N}) + {\rm o}(1).
}

\subsubsection*{Computation of $S^{\rm out}_{2-}$}
Here again, one can expand the $\log$ to deduce
\Eq{}{
 S^{\rm out}_{2-}\underset{N\to\infty}{=}&\sum_{j = \lfloor N\tau_0 \rfloor}^{\lfloor N(\tau_0 + \Delta_N) \rfloor}  \frac{1}{\sqrt{N}}\sqrt{ \frac{2 }{\Delta Q(1)} } \frac{h_1(x(j))}{f_1(x(j))}+{\rm O}(\Delta_N).
}
Explicitly, we have 
\Eq{}{
\frac{h_1(x)}{f_1(x)}=-\frac{1}{3}\left[x(2x^2-3)+(x^2-2)\Phi'(x) \right]-\frac{1}{6}\frac{\partial_r\Delta Q(1)}{\Delta Q(1)}\left[x(2x^2+3)+(x^2+1)\Phi'(x) \right],
}
with
\Eq{eq: asymp h/f}{
\frac{h_1(x)}{f_1(x)}\underset{x\to \infty}{=}-\frac{1}{x}+{\rm O}\left(\frac{1}{x^3}\right).
}
Using the Euler-Maclaurin formula \autoref{prop: EM formula}, and then the change of variable $j\mapsto x(j)$, yields
\Eq{}{
\sum_{j = \lfloor N\tau_0 \rfloor}^{\lfloor N(\tau_0 + \Delta_N) \rfloor}  \frac{1}{\sqrt{N}}\sqrt{ \frac{2 }{\Delta Q(1)} } \frac{h_1(x(j))}{f_1(x(j))}
=\int_0^{\sqrt{\frac{2 }{N\Delta Q(1)}}\xi_N} \left( \frac{h_1(x)}{f_1(x)} + \frac{1}{x + 1} \right) dx - \int_0^{\sqrt{\frac{2 }{N\Delta Q(1)}}\xi_N} \frac{dx}{x + 1}  + {\rm o}(1).
}
Note that we introduced an artificial term $1/(1+x)$ due to the asymptotic \eqref{eq: asymp h/f}. Doing this allows us to extend the integration to $\infty$, on the first integral which is now absolutely convergent. Explicitly,
\Eq{}{
\int_0^{\infty}   \left( \frac{h_1(x)}{f_1(x)} + \frac{1}{x + 1} \right)dx&=\int_0^{\infty}\left( \frac{1}{x + 1} -\frac{1}{3}\left[x(2x^2-3)+(x^2-2)\Phi'(x) \right]\right)dx\\
&\quad -\frac{1}{6}\frac{\partial_r\Delta Q(1)}{\Delta Q(1)}\int_0^{\infty}\left[x(2x^2+3)+(x^2+1)\Phi'(x) \right]dx
}
and 
\Eq{}{
  -\int_0^{\sqrt{\frac{2 }{N\Delta Q(1)}}\xi_N} \frac{dx}{x + 1} \underset{N\to\infty}{=}- \log(\frac{\xi_N}{\sqrt{N}}) -\log \left(\sqrt{\frac{2}{\Delta  Q(1)}}\right)+ {\rm o}(1).
}
Hence, using the definition of $\alpha^{\rm out}$, $\beta^{\rm out}$ \eqref{eq: pos constant out}
\Eq{eq: S2-}{
S^{\rm out}_{2-}\underset{N\to\infty}{=}- \alpha^{\rm out}+\frac{\partial_r\Delta Q(1)}{\Delta Q(1)}\beta^{\rm out}- \log(\frac{\xi_N}{\sqrt{N}})-\log \left(\sqrt{\frac{2}{\Delta  Q(1)}}\right) + {\rm o}(1).
}
\subsubsection*{Computation of $S^{\rm out}$}
Gathering the expansions of $S^{\rm out}_0$ \eqref{eq: S0}, $S^{\rm out}_{1+}$ \eqref{eq: S1+}, $S^{\rm out}_{2+}$ \eqref{eq: S2+}, $S^{\rm out}_{1-}$ \eqref{eq: S1-},$S^{\rm out}_{2-}$ \eqref{eq: S2-}, one finally arrives at
\Eq{4.34}{
S^{\rm out}\underset{N\to\infty}{=}&\,
(- N \eta -  \{N \tau_0\}) (N  q(1) +  \log(N))
-N \eta\left(\log(2\eta)-1 \right) 
    + \xi_N^*\left( \log(\frac{2\xi_N^*}{N})-1\right) \\
   & + \frac{1}{2}\left(\log(2\eta)+\log(\frac{2\xi_N^*}{N})  \right) -\frac{\Delta Q(1)}{4} \left( \frac{N}{\xi_N^*} - \frac{1}{\eta} \right) 
    + \log(\frac{\xi_N^*}{N}) - \log(\eta)\\
&-\gamma^{\rm out}  \sqrt{\Delta Q(1) N}+\frac{1}{4} \log(N)-\frac{1}{4} \log \left(\frac{2 \Delta Q(1)}{\pi }\right)\\
&- \xi_N \left(\log(\frac{2\xi_N}{N})-1\right)+\frac{\Delta Q(1)}{4}\frac{N}{\xi_N}-\frac{1}{2} \log(\frac{2\xi_N}{N})\\
&- \alpha^{\rm out}+\frac{\partial_r\Delta Q(1)}{\Delta Q(1)}\beta^{\rm out}- \log(\frac{\xi_N}{\sqrt{N}})-\log \left(\sqrt{\frac{2}{\Delta  Q(1)}}\right) + {\rm o}(1).
}
This, after simplification using \eqref{eq: def xi_N}, yields the claim
\Eq{eq:2out}{
S^{\rm out}\underset{N\to\infty}{=}&\,
(- N \eta -  \{N \tau_0\}) (N  q(1) +  \log(N))
-\gamma^{\rm out}  \sqrt{\Delta Q(1) N}-N \eta (\log (2 \eta)-1) \\
&-\frac{\log(N)}{4}-\alpha^{\rm out} +\beta^{\rm out}\frac{  \partial_r\Delta Q(1)}{\Delta Q(1)}+\frac{\Delta Q(1)}{4 \eta }-\frac{1}{2} \log(\eta) +\frac{1}{4} \log(\frac{\pi  \Delta Q(1)}{2})+{\rm o}(1).
}
Note how the apparent dependence on
$\Delta_N,\xi_N$ in (\ref{eq: sum2}) in fact cancels out,
 as it should, although a term dependent on $\{N \tau_0\} := N\tau_0 - \lfloor N\tau_0 \rfloor$ remains. We will see that the latter cancels out when added to  the large $N$ expansion of the first sum in (\ref{4.1c}).

\end{proof}

We proceed now with an analogous analysis to that just undertaken for the  first sum in (\ref{4.1c}).

\prop{lem: sum in}{
For $\tau_0 \le 1$, we have the $N\to\infty$ expansion
\begin{multline}
\sum_{j=0}^{\lfloor N\tau_0 \rfloor-1} \log u_j =\,- N^2 I_{\mathbb{S}\cap \mathbb{D}}[\mu_{Q}] - N\log(N) \frac{(1-\eta)}{2} + 
  \{N \tau_0\}(N q(1) + \log(N)) - N \bigg (\frac{1}{2} E_{\mathbb{S}\cap \mathbb{D}}[\mu_{Q}] \\  - \frac{\log(\pi/2)}{2} (1-\eta)\bigg )
-\sqrt{N}\gamma^{\rm in}\sqrt{\Delta Q(1)}+ \frac{\log(2)}{2}+F_{\mathbb{S}\cap \mathbb{D}}\left[Q\right]-\alpha^{\rm in}+\beta^{\rm in}  \frac{\partial_r \Delta Q(1)}{\Delta Q(1)} + {\rm o}(1),
\end{multline}
where the universal positive constants $\alpha^{\rm in}, \beta^{\rm in},\gamma^{\rm in}$ are given in \eqref{eq: pos constant in}.
}
\begin{proof}[Proof of \autoref{lem: sum in}]
As for the  proof of \autoref{lem: sum out}, we start by splitting the sum,
\Eq{}{
S^{\rm in}:=\sum_{j = 0}^{\lfloor N\tau_0 \rfloor -1} \log u_j =&\,\sum_{j = 0 }^{\lfloor N(\tau_0 - \delta_N) \rfloor}\log u_j +\sum_{j = \lfloor N(\tau_0 - \delta_N) \rfloor +1}^{\lfloor N\tau_0 \rfloor -1}\log u_j.
}
The explicit form of the asymptotics of $u_j$ given in \autoref{coro: uj asym4} and \autoref{lem: uj asym5} then gives
\Eq{}{
S^{\rm in}=S^{\rm in}_0+S^{\rm in}_{1-}+S^{\rm in}_{1+}+S^{\rm in}_{2+},
}
where
\Eq{}{
&S^{\rm in}_0:=\sum_{j = 0}^{\lfloor N\tau_0 \rfloor -1} \log(\frac{\sqrt{2\pi}r_\tau}{\sqrt{N V_\tau''(r_\tau)}}e^{-N V_\tau(r_\tau)}),\quad S^{\rm in}_{1-}:=\sum_{j = 0 }^{\lfloor N(\tau_0 - \delta_N)\rfloor}\log(1+\frac{1}{N} \mathcal{B}(r_\tau)),\\
& S^{\rm in}_{1+}:=\sum_{j = \lfloor N(\tau_0 - \delta_N)\rfloor +1}^{\lfloor N\tau_0\rfloor -1} \log(f_2(x)),\quad S^{\rm in}_{2+}:=\sum_{j = \lfloor N(\tau_0 - \delta_N) \rfloor +1}^{\lfloor N\tau_0 \rfloor -1} \log(1+ \sqrt{ \frac{2 }{N\Delta Q(1)} }\frac{h_2(x)}{f_2(x)}),
}
with $x$ defined as in \eqref{eq: def x}.

\subsubsection*{Computation of $S^{\rm in}_{1+}$}
Application of  the Euler-Maclaurin formula \autoref{prop: EM formula} gives
\Eq{eq: EM S1+}{
S^{\rm in}_{1+}\underset{N\to\infty}{=}&\, \sqrt{\frac{N\Delta Q(1)}{2 }} \int_{-\sqrt{\frac{2 }{N\Delta Q(1)}}\zeta_N}^0  \log(f_2(x)) \, dx 
    - \frac{1}{2} \left[ \log(f_2(0)) + \log(f_2\left(-\sqrt{\frac{2 }{N\Delta Q(1)}} \zeta_N\right)) \right] + {\rm o}(1),
}
where $\zeta_N:=\lfloor N(\tau_0 - \delta_N) \rfloor-\lfloor N\tau_0  \rfloor$. The function $x\mapsto\log(f_2(x))$ having the  asymptotic behaviour
\Eq{eq: asymp logf2}{
\log(f_2(x))\underset{x\to -\infty}{\sim} \frac{e^{-x^2}}{2x \sqrt{\pi}},
}
one can extend the domain of integration $[-\sqrt{\frac{2 }{N\Delta Q(1)}} \zeta_N,0]$ to $]-\infty,0]$ with exponentially small error as $N\to\infty$, and also neglect the second boundary term in \eqref{eq: EM S1+}. The first boundary term is given by
\Eq{}{
\log(f_2(0))=-\log(2).
}
Finally, noticing that $\log(f_2(x))=\Phi(x)$ as defined by \eqref{eq: def phi}, we have 
\Eq{eq: EM S1+a}{
S^{\rm in}_{1+}\underset{N\to\infty}{=}&\, \sqrt{\frac{N\Delta Q(1)}{2 }} \int_{-\infty}^0  \Phi(x) \, dx + \frac{\log(2) }{2}  + {\rm o}(1).
}
\subsubsection*{Computation of $S^{\rm in}_{2+}$}
Expanding the $\log$ and applying, once again, the Euler-Maclaurin formula \autoref{prop: EM formula} on $S^{\rm in}_{2+}$ yields
\Eq{eq: EM S2+}{
S^{\rm in}_{2+}\underset{N\to\infty}{=}\,\int_{-\sqrt{\frac{2 }{N\Delta Q(1)}} \zeta_N}^0  \frac{h_2(x)}{f_2(x)} \, dx 
    - \frac{1}{2} \sqrt{ \frac{2 }{N\Delta Q(1)} }\left[ \frac{h_2(0)}{f_2(0)} + \frac{h_2\left(-\sqrt{\frac{2 }{N\Delta Q(1)}} \zeta_N\right)}{f_2\left(-\sqrt{\frac{2 }{N\Delta Q(1)}} \zeta_N\right)} \right] + {\rm o}(1),
}
where 
\Eq{}{
\frac{h_2(x)}{f_2(x)}=-\frac{ \left(x^2-2\right)}{3 }\Phi'(x)-\frac{\partial_r \Delta Q(1)}{\Delta Q(1)}\frac{\left(x^2+1\right) }{6}\Phi'(x).
}
As $\Phi'$ has the asymptotic behaviour
\Eq{}{
\Phi'(x)\underset{x\to-\infty}{\sim}-\frac{e^{-x^2}}{\sqrt{\pi}},
}
one can extend the domain of integration $[-\sqrt{\frac{2 }{N\Delta Q(1)}} \zeta_N,0]$ to $]-\infty,0]$ with exponentially small error as $N\to\infty$ and neglect both boundary terms in \eqref{eq: EM S2+} to get to
\Eq{eq: EM S1+b}{
S^{\rm in}_{2+}\underset{N\to\infty}{=}&\,- \int_{-\infty}^0  \frac{ \left(x^2-2\right)}{3 }\Phi'(x) \, dx-  \frac{\partial_r \Delta Q(1)}{\Delta Q(1)}\int_{-\infty}^0  \frac{\left(x^2+1\right) }{6}\Phi'(x) \, dx + {\rm o}(1).
}
\subsubsection*{Computation of $S^{\rm in}_{1-}$}
One can expand the $\log$ and then extend the sum to $\lfloor N\tau_0\rfloor$ with an error of order ${\rm o}(1)$. This yields
\Eq{}{
S^{\rm in}_{1-}\underset{N\to\infty}{=}\sum_{j = 0 }^{\lfloor N\tau_0\rfloor}\frac{1}{N} \mathcal{B}(r_\tau)+{\rm o}(1).
}
Applying the Euler-Maclaurin formula yields
\Eq{}{
S^{\rm in}_{1-}\underset{N\to\infty}{=}\int_{\mathbb{A}_{r_0,1}} \mathcal{B}(\abs{z}) \, d\mu_Q(z)+{\rm o}(1).
}
The computation of the above integral is given in \cite[Lemma 2.2]{Byun2023} by replacing $r_1$ by $1$.  Explicitly it is
\Eq{}{
\int_{\mathbb{A}_{r_0,1}} \mathcal{B}(\abs{z}) \, d\mu_Q(z)=F_{\mathbb{S}\cap \mathbb{D}}\left[Q\right]-\frac{1}{4}\log(\frac{\Delta Q(1)}{\Delta Q(r_0)})-\frac{1}{3}\log(r_0).
}

\subsubsection*{Computation of $S^{\rm in}_{0}$}
The computation of $S^{\rm in}_{0}$ is also given by replacing $r_1$ by $1$ in \cite[Lemma 2.3 and Lemma 2.4]{Byun2023} and accounting for the facts that we have here $\lfloor N\tau_0 \rfloor$ summands instead of $N$, and our $u_j$ is one half of the corresponding quantity in \cite{Byun2023}. Explicitly, 
\Eq{}{
S^{\rm in}_{0}\underset{N\to\infty}{=}&\,- N^2  I_{\mathbb{S}\cap \mathbb{D}}[\mu_{Q}] - \frac{N}{2} \tau_0\log(N)
+  \{N \tau_0\}(N q(1) + \log(N)) 
+ N \bigg ( \frac{\log(\pi/2)}{2} \tau_0 \\ &- \frac{1}{2} E_{\mathbb{S}\cap \mathbb{D}}[\mu_{Q}] \bigg )
+ \frac{1}{3} \log  r_0 
+ \frac{1}{4} \log\left( \frac{\Delta Q(1)}{\Delta Q(r_0)} \right)
+ {\rm O}\left(N^{-1}\right).
}

\subsubsection*{Computation of $S^{\rm in}$}
Combining, we have

\Eq{}{
S^{\rm in}\underset{N\to\infty}{=}&\,- N^2 I_{\mathbb{S}\cap \mathbb{D}}[\mu_{Q}] - \frac{N}{2} \tau_0\log(N)+ 
  \{N \tau_0\}(N q(1) + \log(N))
+N \Bigg ( \frac{\log(\pi/2)}{2} \tau_0 \\& - \frac{1}{2} E_{\mathbb{S}\cap \mathbb{D}}[\mu_{Q}]  \bigg )
+\sqrt{\frac{N\Delta Q(1)}{2 }} \int_{-\infty}^0  \Phi(x) \, dx + \frac{\log(2)}{2}+F_{\mathbb{S}\cap \mathbb{D}}\left[Q\right]\\
&-\int_{-\infty}^0  \frac{ \left(x^2-2\right)}{3 }\Phi'(x) \, dx-  \frac{\partial_r \Delta Q(1)}{\Delta Q(1)}\int_{-\infty}^0  \frac{\left(x^2+1\right) }{6}\Phi'(x) \, dx + {\rm o}(1),
}
and using that $\tau_0 = 1 - \eta$, the definitions of the constants $\alpha^{\rm in}, \beta^{\rm in},\gamma^{\rm in}$ \eqref{eq: pos constant in} yields the claim.

\end{proof}

\begin{proof}[Proof of \autoref{theo: In/Out Annulus}] ${}$\\
Case $r_1 > 1$. We add together the results of
\autoref{lem: sum out} and \autoref{lem: sum in}. \\
Case $r_1 = 1$. This is the result of \autoref{lem: sum in}.

\end{proof}

\section{Proof In/Out Disk case: $0 = r_0 < 1 \leq r_1$}
As seen in \autoref{S3.3}, in the case $r_0 = 0$ there is a new asymptotic behaviour of $u_j$ relative to the case $r_0 > 0$ for $\tau \le {\rm O}(N^{\epsilon - 1})$, $(0 < \epsilon < 1/5$). It was remarked there that this very same change in asymptotic behaviours occurs in the study \cite{Byun2023} of the partition function $Z_N^{\rm s}[q(r)]$. Since this common alteration is independent of $\tau_0$, we must have
\Eq{Z1}{ Z_N^{\rm s}[q(r)] |_{\rm disk} -
 Z_N^{\rm s}[q(r)] |_{\rm annulus} =
Z_N^{\rm h}[q(r)] |_{\rm disk} -
 Z_N^{\rm h}[q(r)] |_{\rm annulus}. }
 The left hand side of (\ref{Z1}) can be read off from the results of \cite{Byun2023},
 and is equal to 
 $$
 - \frac{1 }{ 12} \log(N)+ \zeta'(-1) + F_{\mathbb D}[Q] - F_{\mathbb{S}\cap \mathbb{D}}[Q],
 $$
 thus implying the stated result.
 
\section{Proof Out Annulus case: $1 \leq  r_0 <r_1$}
We start from the case $1<r_0<r_1$, for which a slight modification of the proof of \autoref{lem: sum out} gives the desired result. We first split the sum as 
\begin{align*}
    S:=\sum_{j=0}^{N-1}\log u_j=S_0+S_1+S_2,
\end{align*}
where
\begin{align*}
    S_0=\sum_{j=0}^{N-1}\log(\frac{\exp(-N V_\tau(1))}{N}),\quad S_1=\sum_{j = 0}^{N-1} \log(\frac{-1}{V_\tau'(1)}),\quad S_2=\sum_{j = 0}^{N-1} \log(1+\frac{1}{N}  \frac{V_\tau'(1) - V_\tau''(1)}{(V_\tau'(1))^2}).
\end{align*}
Use of the Euler-Maclaurin formula gives the analogues of \eqref{eq: S0}, \eqref{eq: S1+} and \eqref{eq: S2+}
\Eq{}{
&S_0=-N^2q(1)-N\log(N),\\
&S_1=N(1-\log(2)-\eta\log\eta+(\eta-1)\log(\eta-1))+\frac{1}{2}\log\frac{\eta}{\eta-1}+{\rm o}(1), \\
&S_2=\log\frac{\eta-1}{\eta}+\frac{\Delta Q(1)}{4\eta(1-\eta)}+{\rm o}(1).
}
Combining these together gives the first asymptotic expansion in \autoref{theo: Out Annulus}. The case $1=r_0<r_1$ is just a special case of \autoref{lem: sum out} with $\tau_0=0$. Therefore, the second formula in \autoref{theo: Out Annulus} is obtained by setting $\eta=1$ in \eqref{eq: sum out}.


\section*{Acknowledgement}
M.A. and S.L. are supported by the International Research Training Group (IRTG) between the University of Melbourne and KU Leuven and Melbourne Research Scholarship of University of Melbourne. P.J.F.~is supported by the Australian Research Council Discovery Project DP250102552. B.S.~is supported by the Australian Research Council Discovery Project DP210102887 and the Shanghai Jiao Tong University Overseas Joint Postdoctoral Fellowship Program. Sung-Soo Byun is to be thanked for providing detailed feedback on an earlier  draft of this work. We also thank Christophe Charlier and Kohei Noda for their comments on our manuscript and verifying our results numerically in relation to their forthcoming work.

\appendix
\section{Appendix}\label{appendix}
Essential to our analysis is that for the class of sums encountered the large $N$ asymptotic expansion can be deduced by applying the Euler-Maclaurin summation formula \cite{DLMF}, which we record here for reference.

\prop{prop: EM formula}{(Euler-Maclaurin formula) Let $f(x)$ be $2k$ times differentiable  on the interval $[p,q]$. Denote by $\{B_{2k} \}$ the even indexed Bernoulli numbers ($B_2=\frac{1 }{ 6}, B_4 = - \frac{1 }{ 30}, \dots$)
\Eq{}{
\sum _{i=p+1}^{q-1}f(i)=\int _{p}^{q}f(x)~{\rm {d}}x-{\frac {f(p)+f(q)}{2}}+\sum _{j=1}^{k}{\frac {B_{2j}}{(2j)!}}\left(f^{(2j-1)}(q)-f^{(2j-1)}(p)\right)+R_{2k}
}
\Eq{}{
\sum _{i=p}^{q}f(i)=\int _{p}^{q}f(x)~{\rm {d}}x+{\frac {f(p)+f(q)}{2}}+\sum _{j=1}^{k}{\frac {B_{2j}}{(2j)!}}\left(f^{(2j-1)}(q)-f^{(2j-1)}(p)\right)+R_{2k},
}
where the remainder $R_{2k}$ satisfies the bound
$|R_{2k}| \le c_{2k}
\int_p^q | f^{(2k)}(x) | \, dx$, with $c_{2k} = 2 \zeta(2k) / (2 \pi)^{2k}$.
}

 \bibliography{main}

\end{document}